\documentstyle[12pt,amsfonts, float, amsmath,amssymb, amsthm, graphicx, multirow, color,leqno,ulem, setspace, caption, subcaption, threeparttable, url, rotating, xr]{article}

\newcommand{\mbf}[1]{\mbox{\boldmath $#1$}}

{\catcode `\@=11 \global\let\AddToReset=\@addtoreset}
\AddToReset{equation}{section}

\AddToReset{Theorem}{section}

\newtheorem{cor}{Corollary}[section]
\newtheorem{lem}{Lemma}[section]
\newtheorem{thm}{Theorem}[section]

\theoremstyle{remark}
\newtheorem{rem}{Remark}[section]

\newcommand{\cC}{{\cal C}}

\newcommand{\cF}{{\cal F}}

\newcommand{\cK}{{\cal K}}
\newcommand{\cL}{{\cal L}}

\newcommand{\cT}{{\cal T}}

\newcommand{\cW}{{\cal W}}

\newcommand{\cZ}{{\cal Z}}

\def\bc{\begin{center}}
\def\bd{\begin{description}}
\def\be{\begin{enumerate}}
\def\ec{\end{center}}
\def\ed{\end{description}}
\def\edt{\end{document}}
\def\ee{\end{enumerate}}
\def\ben{\begin{equation}}
\def\benn{\begin{equation*}}
\def\een{\end{equation}}
\def\eenn{\end{equation*}}
\def\benr{\begin{eqnarray}}
\def\eenr{\end{eqnarray}}
\def\benrr{\begin{eqnarray*}}
\def\eenrr{\end{eqnarray*}}

\def\al{\alpha}
\def\b{\beta}

\def\del{\delta}
\def\edt{\end{document}}

\def\g{\gamma}

\def\iny{\infty}

\def\la{\lambda}
\def\lel{\label}

\def\noi{\noindent}

\def\r{\ref}

\def\ra{\rightarrow}

\def\si{\sigma}

\def\sti{\sum_{i=1}^n}
\def\stj{\sum_{j=1}^n}

\def\t{\tau}

\def\va{\vartheta}
\def\vep{\varepsilon}

\def\vs{\vskip}
\def\wh{\widehat}

\def\wt{\widetilde}

\def\R{{\mathbb R}}
\def\Z{{\mathbb Z}}
\def\z{\zeta}

\begin{document}

\bc
{\large \bf Generalized Minimum Distance Estimators in Linear \\Regression with Dependent Errors}\\[.2cm]
{\large Jiwoong Kim\\University of Notre Dame}
\ec
\vs .75cm

\begin{abstract}
\vs .1cm
This paper discusses minimum distance estimation method in the linear regression model with dependent errors which are strongly mixing. The regression parameters are estimated through the minimum distance estimation method, and  asymptotic distributional properties of the estimators are discussed. A simulation study compares the performance of the minimum distance estimator with other well celebrated estimator. This simulation study shows the superiority of the minimum distance estimator over another estimator. \texttt{KoulMde} (\texttt{R} package) which was used for the simulation study is available online. See section \r{sim} for the detail.
\end{abstract}
\vs 0.75cm
\begin{minipage}{1\textwidth}
$Keywords$: Dependent errors; Linear regression; Minimum distance estimation; Strongly mixing
\end{minipage}

\section{Introduction}
\noindent
Consider the linear regression model
\benr\label{eq:Model1}
y_{i} =  \mbf{x}_i'\mbf{\b} + \vep_i,
\eenr
where $E\vep_i\equiv 0$, $\mbf{x}_i=(1, x_{i2},...,x_{ip})'\in \R^{p}$ with $x_{ij},
j=2,\cdots, p,\, i= 1,\cdots, n$ being non random design variables, and where
$\mbf{\b} = (\beta_{1},...,\beta_{p})'\in \R^{p}$ is the parameter vector of interest. The methodology where the estimators are obtained by minimizing some dispersions or pseudo distances between the data and the underlying model is referred to as the minimum distance (m.d.) estimation method. In this paper we estimate regression parameter vector $\mbf{\b}$ by the m.d.\,\,estimation method when the collection of $\vep_i$ in the model (\r{eq:Model1}) is a dependent process.

Let $T_{1},...,T_{n}$ be independent identically distributed (i.i.d.) random variables (r.v.'s) with distribution function (d.f.) $G_{\vartheta}$ where $\vartheta$ is unknown. The classical m.d.\,\,estimator of $\vartheta$ is obtained by minimizing following Cram$\grave{\textrm{e}}$r-von Mises (CvM) type $L_2$-distance
\ben\lel{eq:cvm}
\int \big\{\,G_{n}(y) - G_{\vartheta}(y)    \big\}^2 dH(y)
\een
where $G_{n}$ is an empirical d.f. of $T_{i}$'s and $H$ is a integrating measure. There are multiple reasons as to why CvM type distance is preferred, including the asymptotic normality of the corresponding m.d.\,\,estimator; see, e.g., Parr and Schucany (1980), Parr (1981) and Millar (1981). Many researchers have tried various $H$'s to obtain the m.d.\,\,estimators. Anderson and Darling (1952) proposed Anderson-Darling estimator obtained by using $dG_{\vartheta}/\{G_{\vartheta}(1-G_{\vartheta})\}$ for $dH$. Another important example includes $H(y)\equiv y$, giving a rise to Hodges - Lehmann type estimators. If $G_{n}$ and $G_{\vartheta}$ of the integrand are replaced with kernel density estimator and assumed density function of $T_{i}$, the Hellinger distance estimators will be obtained; see Beran (1977).

Departing from one sample setup, Koul and DeWet (1983) extended the domain of the application of the m.d.\,\,estimation to the regression setup. On the assumption that $\vep_i$'s are i.i.d.\,\,r.v.'s with a known d.f.\,\,$F$, they proposed a class of the m.d.\,\,estimators by minimizing $L_2$-distances between a weighted empirical d.f.\,\,and the error d.f.\,\,$F$. Koul (2002) extended this methodology to the case where error distribution is unknown but symmetric around zero. Furthermore, it was shown therein that when the regression model has independent non-Gaussian errors the m.d.\,\,estimators of the regression parameters --- obtained by minimizing $L_{2}$-distance with various integrating measures --- have the least asymptotic variance among other estimators including Wilcoxon rank, the least absolute deviation (LAD), the ordinary least squares (OLS) and normal scores estimators of $\mbf{\b}$: e.g.\,\,the m.d.\,\,estimators obtained with a degenerate integrating measure display the least asymptotic variance when errors are independent Laplace r.v's.

However, the efficiency of the m.d.\,\,estimators depends on the assumption that errors are independent; with the errors being dependent, the m.d.\,\,estimation method will be less efficient than other estimators. Examples of the more efficient methods include the generalized least squares (GLS); GLS is nothing but regression of transformed $y_{i}$ on transformed $\mbf{x}_{i}'$. The most prominent advantage of using the GLS method is ``decorrelation" of errors as a result of the transformation. Motivated by efficiency of m.d.\,\,estimators --- which was demonstrated in the case of independent non-Gaussian errors --- and the desirable property of the GLS (decorrelation of the dependent errors), the author proposes generalized m.d.\,\,estimation method which is a mixture of the m.d.\,\,and the GLS methods: the m.d.\,\,estimation will be applied to the transformed variables. ``Generalized" means the domain of the application of the m.d.\,\,method covers the case of dependent errors; to some extent, the main result of this paper generalizes the work of Koul (2002). As the efficiency of the m.d.\,\,method is demonstrated in the case of independent errors, the main goal of this paper is to show that the generalized m.d.\,\,estimation method is still competitive when the linear regression model has dependent errors; indeed, the simulation study empirically shows that the main goal is achieved.

The rest of this article is organized as follows. In the next section, characteristics of dependent errors used through this paper is studied. Also, the CvM type distance and various processes --- which we need in order to obtain the estimators of $\mbf{\b}$ --- will be introduced. Section \r{asnb} describes the asymptotic distributions and some optimal properties of the estimators. Findings of a finite sample simulations are described in Section \r{sim}. All the proofs are deferred until Appendix. In the remainder of the paper, an Italic and boldfaced variable denotes a vector while a non-Italic and boldfaced variable denotes a matrix. An identity matrix will carry a suffix showing its dimension: e.g. $\textbf{I}_{n\times n}$ denotes a $n\times n$ identity matrix. For a function $f:\R\ra \R$, let $|f|_{H}^{2}$ denote $\int\,f^{2}(y)\,dH(y)$. For a real vector $\mbf{u}\in\R^{p}$, $\|\mbf{u}\|$ denotes Euclidean norm. For any r.v.\,\,$Y$, $\|Y\|_{p}$ denotes $(E |Y|^{p})^{1/p}$. For a real matrix $\textbf{W}$ and $y\in\R$, $\textbf{W}(y)$ means that its entries are functions of $y$.

\section{Strongly mixing process \& CvM type distance}\lel{esr}
Let $\cF_{m}^{l}$ be the $\sigma$-field generated by $\vep_{m},\vep_{m+1},...,\vep_{l},\, m\le l$. The sequence $\{\vep_{j}, \,j\in \Z\}$ is said to satisfy the strongly mixing condition if
\benn
\alpha(k):=\sup\left\{ |P(A\cap B)-P(A)P(B)|:\,A\in \cF_{-\iny}^{0}, \,B\in \cF_{k}^{\iny}
\right\}\ra 0,
\eenn
as $k\ra\iny$. $\al$ is referred to as mixing number. Chanda (1974), Gorodetskii (1977), Koul (1977), and Withers (1979) investigated the decay rate of the mixing number. Having roots in their works, Section \r{asnb} defines the decay rate assumed in this paper; see, e.g., the assumption (a.8). Hereinafter the errors $\vep_i$'s are assumed to be strongly mixing with mixing number $\al$. In addition, $\vep_i$ is assumed to be stationary and symmetric around zero.

Next, we introduce the basic processes and the distance which are required to obtain desired result. Recall the model (\r{eq:Model1}). Let $\textbf{X}$ denote the $n\times p$ design matrix whose $i$th row vector is $\mbf{x}_{i}'$. Then the model (\r{eq:Model1}) can be expressed as
\benn
\mbf{y} = \textbf{X}\mbf{\b}+\mbf{\vep},
\eenn
where $\mbf{y}=(y_{1},y_{2},...,y_{n})'\in\R^{n}$ and $\mbf{\vep}=(\vep_{1},\vep_{2},...,\vep_{n})'\in\R^{n}$
Let $\textbf{Q}$ be any $n\times n$ real matrix so that the inverse of $\textbf{Q}^2$ is a positive definite symmetric matrix. Note that the diagonalization of positive definite symmetric matrix guarantees the existence of $\textbf{Q}$ which is also a symmetric matrix. Let $\mbf{q}_{i}'=(q_{i1},...,q_{in})$ for $1\leq i\leq n$ denote the $i$th row vector of $\textbf{Q}$. Define transformed variables
\benn
\wt{y}_{i} = \mbf{q}_{i}'\mbf{y},\quad \wt{\mbf{x}}_{2}'=\mbf{q}_{i}'\textbf{X},\quad \wt{\vep}_{i} = \mbf{q}_{i}'\mbf{\vep},\quad 1\le i\le n.
\eenn
As in the GLS method, $\textbf{Q}$ obtained from covariance matrix of $\mbf{\vep}$ transforms dependent errors into uncorrelated ones, i.e., ``decorrelates" the errors. However, the GLS obtains $\textbf{Q}$ in a slightly different manner. Instead of using $\textbf{Q}^{2}$, the GLS equates $\textbf{Q}'\textbf{Q}$ to the inverse of the covariance matrix, i.e., the GLS uses Cholesky decomposition. The empirical result in Section \r{sim} describes that $\textbf{Q}$ from the diagonalization yields better estimators. Here we propose the class of the generalized m.d.\,\,estimators of the regression parameter upon varying $\textbf{Q}$. We impose Noether (1949) condition on $\textbf{QX}$. Now let $\textbf{A}=(\textbf{X}'\textbf{Q}^{2}\textbf{X})^{-1/2}$ and $\mbf{a}_{j}$ denote $j$th column of $\textbf{A}$. Let $\textbf{D}=((d_{ik}))$, $1\leq i\leq n$, $1\leq k\leq p$, be an $n\times p$ matrix of real numbers and $\mbf{d}_j$ denote $j$th column of $\textbf{D}$. As stated in Koul (2002, p.60), if $\textbf{D}=\textbf{QXA}$ ( i.e., $d_{ik} = \mbf{q}_{i}'\textbf{X}\mbf{a}_{k}$), then under Noether condition,
\ben\lel{eq:dik2}
\sum_{i}^{n}d_{ik}^{2}= 1,\,\,\, \max_{1\leq i\leq n}d_{ik}^{2}= o(1)\quad\quad \textrm{ for all }1\leq k\leq p.
\vspace{-.4cm}
\een
Next, define CvM type distance from which the generalized m.d.\,\,estimator are obtained. Let $f_{i}$ and $F_{i}$ denote the density function and the d.f.\,\,of $\wt{\vep}_{i}$, respectively. Analogue of (\r{eq:cvm}) --- with $G_{n}$ and $G_{\vartheta}$ being replaced by empirical d.f.\,\,of $\wt{\vep}_{i}$ and $F_{i}$ --- will be a reasonable candidate. However, the d.f.\,\,$F_{i}$ is rarely known. Since the original regression error $\vep_{i}$'s are assumed to be symmetric, the transformed error $\wt{\vep}_{i}$'s are also symmetric; therefore we introduce, as in Koul (2002; Definition 5.3.1),
\benrr\lel{eq:intro1_1}
U_{k}(y,\mbf{b};\, \textbf{Q}) &:=& \sti d_{ik}\Big\{\,I\big( \mbf{q}_{i}'\mbf{y}-\mbf{q}_{i}'\textbf{X}\mbf{b} \leq y\big) -
I\big(-\mbf{q}_{i}'\mbf{y}+\mbf{q}_{i}'\textbf{X}\mbf{b} < y\big)     \,\Big\}, \\
\mbf{U}(y,\mbf{b};\, \textbf{Q})&:=& (U_{1}(y,\mbf{b};\, \textbf{Q}),...,U_{p}(y,\mbf{b};\, \textbf{Q}))', \quad y\in \R,\\
\cL(\mbf{b};\, \textbf{Q}) &:=& \int \|\mbf{U}(y,\mbf{b};\, \textbf{Q})\|^{2} \, d H(y), \quad\quad\quad \mbf{b}\in \R^{p},\\
     &=& \sum_{k=1}^{p}\int\left[\sti d_{ik}\Big\{\,I\big( \mbf{q}_{i}'\mbf{y}-\mbf{q}_{i}'\textbf{X}\mbf{b} \leq y\big) -
I\big(-\mbf{q}_{i}'\mbf{y}+\mbf{q}_{i}'\textbf{X}\mbf{b} < y\big)     \,\Big\}\right]^{2},
\eenrr
where $I(\cdot)$ is an indicator function, and $H$ is a $\sigma-$finite measure on $\R$ and symmetric around 0, i.e., $dH(-x)=-dH(x), \, x\in \R$. Subsequently, define $\wh{\mbf{\b}}$ as
\benn
\cL(\wh{\mbf{\b}};\, \textbf{Q})=\inf_{\substack{b\in\R^{p}}}\cL(\mbf{b};\, \textbf{Q}).
\eenn
Next, define
\benn
                \wt{\textbf{Q}}_{i}:=\left[
                   \begin{array}{c}
                     \mbf{q}_{i}'\\
                     \mbf{0}' \\
                     \vdots \\
                     \mbf{0}' \\
                   \end{array}
                 \right],\quad \wt{\textbf{Q}}:=\left[
                   \begin{array}{c}
                     \wt{\textbf{Q}}_{1}\\
                     \wt{\textbf{Q}}_{2} \\
                     \vdots \\
                     \wt{\textbf{Q}}_{n} \\
                   \end{array}
                 \right],
\eenn
where $\mbf{q}_{i}'$ is the $i$th row vector of $\textbf{Q}$ and $\mbf{0}=(0,...,0)'\in\R^{n}$; observe that $\wt{\textbf{Q}}_{i}$ and $\wt{\textbf{Q}}$ are $n\times n$ and $n^2\times n$ matrices, respectively. Define a $n\times n^{2}$ matrix $\textbf{I}_{f}(y)$ so that its $(i,j)$th entry is $f_{i}(y)I(j=n(i-1)+1)$: e.g., $(k,k(k-1)+1)$th entry is $f_{k}(y)$ for all $1\le k\le n$ and all other entries are zeros. Finally, define following matrices:
\ben\lel{eq:Sig}
\mbf{\Sigma}_{\textbf{D}}:=\int\,\textbf{I}_{f}'(y)\textbf{D}\textbf{D}'\textbf{I}_{f}(y)\,dH(y),\quad\mbf{\Sigma}:=\textbf{A}\textbf{X}'\wt{\textbf{Q}}'\mbf{\Sigma}_{\textbf{D}}\wt{\textbf{Q}}\textbf{X}\textbf{A},
\een
which are needed for the asymptotic properties of $\wh{\mbf{\b}}$. Let $f_{ij}^{H}:=\int\,f_{i}f_{j}\, dH$ and $d_{ij}^{*}:=\sum_{k=1}^{p}d_{ik}d_{jk}$. Note
that
\benn
\mbf{\Sigma}=\textbf{A}\textbf{X}' \left[  \sti\stj d_{ij}^{*}f_{ij}^{H}\mbf{q}_{i}\mbf{q}_{j}' \right] \textbf{X}\textbf{A}.
\eenn

\section{Asymptotic distribution of $\wh{\mbf{\b}}$}\lel{asnb}
In this section we investigate the asymptotic distribution of $\wh{\mbf{\b}}$ under the current setup. Note that minimizing $\cL(\cdot;\,\textbf{Q})$ does not have the closed form solutions; only numerical solutions can be tried, and hence it would be impracticable to derive asymptotic distribution of $\wh{\mbf{\b}}$. To redress this issue, define for $\mbf{b}\in \R^{p}$
\benn
{\cL}^{*}(\mbf{b};\,\textbf{Q}) = \int\,\big\|\mbf{U}(y,\mbf{\b};\,\textbf{Q}) + 2\mbf{\Sigma}_{\textbf{DA}}(y)\textbf{A}^{-1}(\mbf{b}-\mbf{\b}) \big\|^{2}\,d H(y),
\eenn
where $\mbf{\Sigma}_{\textbf{DA}}(y):=\textbf{D}'\textbf{I}_{f}(y)\wt{\textbf{Q}}\textbf{X}\textbf{A}$ is a $p\times p$ matrix. Next, define
\benn
{\cL}^{*}(\widetilde{\mbf{\b}};\,\textbf{Q})=\inf_{\substack{b\in\R^{p}}}{\cL}^{*}(\mbf{b};\,\textbf{Q}).
\eenn
Unlike $\cL(\cdot;\,\textbf{Q})$, minimizing $\cL^{*}(\cdot;\,\textbf{Q})$ has the closed form solution. Therefore, it is not unreasonable to approximate the asymptotic distribution of $\wh{\mbf{\b}}$ by one of $\widetilde{\mbf{\b}}$ if $\cL(\cdot;\,\textbf{Q})$ can be approximated by $\cL^{*}(\cdot;\,\textbf{Q})$. This idea is plausible under certain conditions which are called \textit{uniformly locally asymptotically quadratic}; see Koul (2002, p.159) for the detail. Under these conditions, it was shown that difference between $\wh{\mbf{\b}}$ and $\widetilde{\mbf{\b}}$ converges to zero in probability; see theorem 5.4.1. The basic method of deriving the asymptotic properties of $\wh{\mbf{\b}}$ is similar to that of sections 5.4, 5.5 of Koul (2002). This method amounts to showing that $\cL(\mbf{\b} + \textbf{A}\mbf{u};\,\textbf{Q})$ is uniformly locally asymptotically quadratic in $\mbf{u}$ belonging to a bounded set and  $\|\textbf{A}^{-1}(\wh{\mbf{\b}}-\mbf{\b})\|=O_p(1).$ To achieve these goals we need the following assumptions which in turn have roots in section 5.5 of Koul (2002).
\begin{description}
   \item[(a.1)] The matrix $\textbf{X}'\textbf{Q}^{2}\textbf{X}$ is nonsingular and, with $\textbf{A}=(\textbf{X}'\textbf{Q}^{2}\textbf{X})^{-1/2}$, satisfies
   \benn
   \limsup_{n\rightarrow\iny}\, n\max_{1\leq j\leq p}\|\mbf{d}_{j}\|^2<\iny.
   \eenn
   \item[(a.2)] The integrating measure $H$ is $\sigma-$finite and symmetric around 0, and
   \benn
   \int_{0}^{\iny} (1-F_{i})^{1/2}dH <\iny,\quad 1\le i\le n.
   \eenn
   \item[(a.3)] For any real sequences $\{a_{n}\}$, $\{b_{n}\}$, $b_{n}-a_{n}\rightarrow 0$,
   \benn
   \limsup_{n\rightarrow\iny} \int_{a_{n}}^{b_{n}}\int f_{i}(y+x)dH(y)dx=0,\quad 1\le i\le n.
   \eenn

   \item[(a.4)] For $a\in \R$, define $a^{+}:=\max(a,0)$, $a^{-}:=a^{+}-a$. Let $\theta_{i}:=\|\mbf{q}_{i}'\textbf{XA}\|$. For all $\mbf{u}\in\R^{p}$, $\|\mbf{u}\|\leq b$, for all $\delta>0$, and for all $1\leq k\leq p$,
   \benn
   \limsup_{n\rightarrow\iny} \int\Big[ \sum_{i=1}^{n}d_{ik}^{\pm}\big\{ F_{i}(y+\mbf{q}_{i}'\textbf{XA}\mbf{u} + \delta\theta_{i})
    -F_{i}(y+\mbf{q}_{i}'\textbf{XA}\mbf{u}-\delta\theta_{i})\big\}\Big]^{2}d H(y)\leq c\delta^{2},
   \eenn
   where $c$ does not depend on $\mbf{u}$ and $\delta$.
   \item[(a.5)] For each $u\in \R^{p}$ and all $1\leq k\leq p$,
   \benn
   \int \Big[ \sum_{i=1}^{n}d_{ik}\big\{ F_{i}(y+\mbf{q}_{i}'\textbf{XA}\mbf{u}) -F_{i}(y) - \mbf{q}_{i}'\textbf{XA}\mbf{u}f_{i}(y)\big\}\Big]^{2}d H(y)=o(1).
   \eenn
  \item[(a.6)] $F_{i}$ has a continuous density $f_{i}$ with respect to the Lebesgue measure on $(\mathbb{R},\,\mathcal{B})$ for $i=1,2,...,n$.
  \item[(a.7)] $0<\int_{0}^{\iny} f_{i}^{r}d H<\iny$, for $r=1/2,1,2$ and $i=1,2,...,n$.
  \item[(a.8)] The $\{\vep_{i}\}$ in the model (\r{eq:Model1}) is strongly mixing with mixing number $\alpha(\cdot)$ satisfying
      \benn
      \limsup_{n\ra \iny} \sum_{k=1}^{n-1}k^{2}\alpha(k)<\iny.
      \eenn

\end{description}

\begin{rem}\lel{rem:1}
Note that (a.1) implies Noether condition and (a.2) implies $\int_{0}^{\iny} (1-F_{i}) dH < \iny$.
From Corollary 5.6.3 of Koul (2002), we note that in the case of i.i.d.\,\,errors, the asymptotic
normality of $\wh{\mbf{\b}}$ was established under the weaker conditions: Noether condition and $\int_{0}^{\iny}
(1-F_{i})dH < \iny$.  The dependence of the errors now forces us to assume
two stronger conditions (a.1) and (a.2).
\end{rem}

\begin{rem}\lel{rem:1_1}
Here we discuss examples of $H$ and $F$ that satisfy (A.2). Clearly it is satisfied by any finite measure $H$. Next consider the $\si$-finite measure $H$ given by $d H\equiv \{F_{i}(1-F_{i})\}^{-1}dF_{i}$, $F$ a continuous d.f. symmetric around zero. Then $F_{i}(0) = 1/2$ and
\benn
\int_{0}^{\iny} (1-F_{i})^{1/2}d H = \int_{0}^{\iny} \frac{ (1-F_{i})^{1/2}}{ F_{i}(1-F_{i}) }  dF_{i}
 \leq 2 \int_{1/2}^{1} (1-u)^{-1/2} du < \iny.
\eenn
Another useful example of a $\si$-finite measure $H$ is given by $H(y)\equiv y$. For this measure, (a.2)
is satisfied by many symmetric error d.f.s including normal, logistic, and Laplace. For example, for normal d.f., we do not have a closed form of the integral, but by using the well celebrated tail bound for normal distribution --- see e.g., Theorem 1.4 of Durrett (2005) --- we obtain
\benn
\int_{0}^{\iny} \{1-F_{i}(y)\}^{1/2}dy \leq  (2\pi)^{-1/2}\int_{0}^{\iny} y^{-1/2}\exp(-y^{2}/4) dy = (2/\pi)^{1/2}\,\Gamma(1/4).
\eenn
 Recall from Koul (2002) that the $\wh{\mbf{\b}}$ corresponding to $H(y)\equiv y$ is the extensions of the
 one sample Hodges-Lehmann estimator of the location parameter to the above regression model. 
\end{rem}

\begin{rem}\lel{rem:2} 
Consider condition (a.7).
 If $f_{i}$'s are bounded then $\int f_{i}^{1/2} dH<\iny$ implies the other two conditions in (a.7) for any
 $\si$-finite measure $H$. For $H(y)\equiv y$, $\int f_{i}^{1/2}(y) dy<\iny$ when $f_{i}$'s are normal,
logistic or Laplace densities. In particular, when $d H= \{F_{i}(1-F_{i})\}^{-1}dF_{i}$  and  $F_{i}$'s are
logistic d.f.'s, so that $dH(y)\equiv dy$, this condition is also satisfied. 
\end{rem}

We are ready to state the needed results. The first theorem establishes the needed uniformly locally asymptotically quadraticity while the corollary shows the boundedness of a suitably standardized $\wh{\mbf{\b}}$. Theorem 3.1 and Corollary 3.1 are counterparts of conditions (A$\tilde{1}$) and (A5) in theorem 5.4.1 of Koul (2002), respectively. Note that condition (A4) in theorem 5.4.1 is met by (\r{eq:A4}) in the Appendix; condition (A6) in theorem 5.4.1 is trivial.
\begin{thm}\lel{thm:1}
\textit{Let $\{y_{i},\,1\leq i\leq n\}$ be in the model (\r{eq:Model1}). Assume that (a.1)-(a.8) hold.
Then, for any $0<c<\iny$,}
\ben\lel{eq:thm1}
E\sup_{\|A^{-1}(b-\b)\|\le c} \|\cL(\mbf{b};\,\textbf{Q})-\cL^{*}(\mbf{b};\,\textbf{Q}) \|=o(1).
\een
\end{thm}
\noindent
\textbf{Proof.} See Appendix. \qed
\begin{cor}\lel{cor:1}
\textit{Suppose that the assumptions of Theorem \r{thm:1} hold. Then for any $\epsilon>0$, $0<M<\iny$ there exists an $N_{\vep}$, and $0<c_{\epsilon}<\iny$ such that}
\ben\lel{eq:cor1}
P\left(\inf_{\|A^{-1}(b-\b)\|\ge c_{\epsilon}}\cL(\mbf{b};\,\textbf{Q})\geq M \right)\geq 1-\epsilon,\quad\quad \forall\, n\geq N_{\epsilon}.
\een
\end{cor}
\noindent
\textbf{Proof.} See Appendix. \qed
\begin{thm}\lel{thm:2}
\textit{Under the assumptions of Theorem \r{thm:1},}
\benn 
\textbf{A}^{-1}(\widehat{\mbf{\b}} - \mbf{\b})= -\frac{1}{2} \mbf{\Sigma}^{-1}
\textbf{A}\textbf{X}'\wt{\textbf{Q}}' \int\,\textbf{I}'_{f}(y)\textbf{D} \mbf{U}(y,\mbf{\b})\,d H(y)+o_{p}(1),
\eenn
\textit{where $\mbf{\Sigma}$ is as in (\r{eq:Sig}).}
\end{thm}
\noindent
\textbf{Proof}. Note that the first term in the right-hand side is nothing but $\textbf{A}^{-1}(\wt{\mbf{\b}} - \mbf{\b})$. Therefore, the proof follows from Theorem \r{thm:1} and Corollary \r{cor:1}, as in i.i.d.\,\,case illustrated in the theorem 5.4.1 of Koul (2002). \qed

Next, define
\benrr\lel{eq:KFH}
&&\psi_{i}(x):= \int_{-\iny}^{-x}\,f_{i}(y)\,d H(y) - \int_{-\iny}^{x}\,f_{i}(y)\,d H(y),\\
&&\textbf{Z}_{n} := \int\,\textbf{I}'_{f}(y)\textbf{D} \mbf{U}(y,\mbf{\b})\,d H(y).
\eenrr
Symmetry of the $F_{i}$ around 0 yields $E\psi_{i}(\wt{\vep}_{j})=0$ for $1\le i,j\le n$.
Let $\mbf{\Sigma}_{\textbf{ZZ}}$ denote covariance matrix of $\textbf{A}\textbf{X}'\wt{\textbf{Q}}'\textbf{Z}_{n}$. Define a $n\times n$ matrix $\mbf{\Sigma}_{\psi}$ and write
$\mbf{\Sigma}_{\psi}=((\g_{ij}))$, $1\le i\le p$, $1\le j\le p$ where
\benn
\g_{ij}=\sum_{l=1}^{n}\sum_{h=1}^{n}d_{il}^{*}d_{jh}^{*}E\left[ \psi_{i}( \mbf{q}_{l}'\mbf{\vep})\psi_{j}(\mbf{q}_{h}'\mbf{\vep}) \right].
\eenn
Observe that
\benrr
&&\mbf{\Sigma}_{\textbf{ZZ}}= E(\textbf{A}\textbf{X}'\wt{\textbf{Q}}'\textbf{Z}_{n}\textbf{Z}_{n}'\wt{\textbf{Q}}\textbf{X}\textbf{A}) =  \textbf{A}\textbf{X}'\textbf{Q}'\mbf{\Sigma}_{\psi}\textbf{Q}\textbf{X}\textbf{A}.
\eenrr
Now, we are ready to state the asymptotic distribution of $\wh{\mbf{\b}}$.
\begin{lem}\lel{lem:MR}
\textit{Assume $\mbf{\Sigma_{\textbf{ZZ}}}$ is positive definite for all $n\ge p$. In addition, assume that}
\benn
\sup_{u\in \R^{p}, \|u\|=1} \mbf{u}'\mbf{\Sigma_{\textbf{ZZ}}}^{-1}\mbf{u}=O(1).
\eenn
\textit{Then}
\benn
\boldsymbol{\Sigma}_{\textbf{ZZ}}^{-1/2} \textbf{\textrm{A}}\textbf{X}'\wt{\textbf{Q}}'\textbf{Z}_{n}\to_{d} N(\mbf{0},\textbf{I}_{p\times p}),
\eenn
\textit{where $\mbf{0}=(0,...,0)'\in\R^{p}$ and }$\textbf{I}_{p\times p}$ \textit{is the $p\times p$ identity matrix.}
\end{lem}
\noindent
\textbf{Proof.} To prove the claim, it suffices to show that for any $\mbf{\la}\in\R^{p}$, $\mbf{\la}'\mbf{\Sigma_{\textbf{ZZ}}}^{-1/2} \textbf{A}\textbf{X}'\wt{\textbf{Q}}'\textbf{Z}_{n}$ is asymptotically normally distributed. Note that
\benn
\mbf{\la}'\mbf{\Sigma_{\textbf{ZZ}}}^{-1/2} \textbf{A}\textbf{X}'\wt{\textbf{Q}}'\textbf{Z}_{n} = \sti\left[ \left( \mbf{\la}'\mbf{\Sigma_{\textbf{ZZ}}}^{-1/2} \textbf{A}\textbf{X}'\mbf{q}_{i} \right) \stj d_{ij}^{*}\psi_{i}(\mbf{q}_{j}'\mbf{\vep}) \right],
\eenn
which is the sum as in the theorem 3.1 from Mehra and Rao (1975) with $c_{ni}=\mbf{\la}'\mbf{\Sigma_{\textbf{ZZ}}}^{-1/2} \textbf{A}\textbf{X}'\mbf{q}_{i}$ and $\xi_{ni} = \stj d_{ij}^{*}\psi_{i}(\mbf{q}_{j}'\mbf{\vep})$. Note that
\benn
\t_{c}^{2}:= \sti c_{ni}^{2}= \mbf{\la}'\mbf{\Sigma}_{\textbf{ZZ}}^{-1}\mbf{\la},\,\,\si_{n}^{2}:= E\left\{ \sti\left[ \left( \mbf{\la}'\mbf{\Sigma_{\textbf{ZZ}}}^{-1/2} \textbf{A}\textbf{X}'\mbf{q}_{i} \right) \stj d_{ij}^{*}\psi_{i}(\mbf{q}_{j}'\mbf{\vep}) \right] \right\}^2 = \|\mbf{\la}\|^{2}.
\eenn
Also, observe that
\benn
\max_{1\leq i\leq n} c_{ni}^{2}/ \t_{c}^{2}
\leq \max_{1\leq i\leq n} \frac{  \|\mbf{\la}'\mbf{\Sigma}_{\textbf{ZZ}}^{-1}\|^{2} \| \textbf{A}\textbf{X}'\mbf{q}_{i} \|^{2}}{\mbf{\la}'\mbf{\Sigma}_{\textbf{ZZ}}^{-1}\mbf{\la}}=\max_{1\leq i\leq n} \|\textbf{A}\textbf{X}'\mbf{q}_{i}\|^2 \ra 0,
\eenn
by assumption (A.1). Finally, we obtain
\benn
\liminf_{n\ra \iny}\sigma_{n}^{2}/\tau_{c}^{2} \geq \|\mbf{\la}\|^{2}/(\limsup \mbf{\la}'\mbf{\Sigma}_{\textbf{ZZ}}^{-1}\mbf{\la}) > 0,
\eenn
by the assumption that the terms in the denominator is $O(1)$. Hence, the desired result follows from the theorem 3.1 of Mehra and Rao (1975). \qed

\begin{cor}\lel{cor:MR2}
\textit{In addition to the assumptions of Theorem \r{thm:1}, let the assumption of Lemma \r{lem:MR} hold. Then}
\benn
\mbf{\Sigma}_{\textbf{ZZ}}^{-1/2} \mbf{\Sigma}\textbf{A}^{-1}(\widehat{\mbf{\b}} - \mbf{\b})\ra_{d} 2^{-1}\,N(\mbf{0}, \textbf{I}_{p\times p}).
\eenn
\end{cor}
\noindent
\textbf{Proof}. Claim follows from Lemma \r{lem:MR} upon noting that
\benn
 \textbf{Z}_{n} = \int\,\textbf{I}'_{f}(y)\textbf{D} \mbf{U}(y,\mbf{\b})\,d H(y).
\eenn
\qed

\begin{rem}\lel{Remark:asympVar}
Let $Asym(\wh{\mbf{\b}})$ denote the asymptotic variance of $\wh{\mbf{\b}}$. Then we have
\benrr
Asym(\wh{\mbf{\b}}) &=& 4^{-1}\textbf{A}\mbf{\Sigma}^{-1} \mbf{\Sigma}_{\textbf{ZZ}} \mbf{\Sigma}^{-1}\textbf{A}\\
&=& 4^{-1}\textbf{A} (\textbf{A}\textbf{X}'\wt{\textbf{Q}}'\mbf{\Sigma}_{\textbf{D}}\wt{\textbf{Q}}\textbf{X}\textbf{A})^{-1} (\textbf{A}\textbf{X}'\textbf{Q}'\mbf{\Sigma}_{\psi}\textbf{Q}\textbf{X}\textbf{A}) (\textbf{A}\textbf{X}'\wt{\textbf{Q}}'\mbf{\Sigma}_{\textbf{D}}\wt{\textbf{Q}}\textbf{X}\textbf{A})^{-1}\textbf{A}
\eenrr
Observe that if all the transformed errors have the same distribution, i.e., $f_{1}=f_{2}=\cdots=f_{n}$, we have
\benn
\textbf{A}\textbf{X}'\wt{\textbf{Q}}'\mbf{\Sigma}_{\textbf{D}}\wt{\textbf{Q}}\textbf{X}\textbf{A}=(|f_{1}|_{H}^{2})^{-1}\textbf{I}_{p\times p}.
\eenn
Therefore, $Asym(\wh{\mbf{\b}})$ will be simplified as
\benn
(2|f_{1}|_{H}^{2})^{-2}\textbf{A}(\textbf{A}\textbf{X}'\textbf{Q}'\mbf{\Sigma}_{\psi}\textbf{Q}\textbf{X}\textbf{A})\textbf{A}.
\eenn
Moreover, if all the transformed errors are uncorrelated as a result of the transformation, $Asym(\wh{\mbf{\b}})$ can be simplified further as
\benn
\tau(2|f_{1}|_{H}^{2})^{-2} (\textbf{X}'\textbf{Q}^{2}\textbf{X})^{-1},
\eenn
where $\tau= Var(\psi_{1}(\mbf{q}_{1}'\mbf{\vep}))$.
\end{rem}


\section{Simulation studies}\lel{sim}
\noi
In this section the performance of the generalized m.d.\,\,estimator is compared with one of the GLS estimators. Let $\mbf{\Omega}:=E(\mbf{\vep}\mbf{\vep}')$ and $\wh{\mbf{\Omega}}$ denote covariance matrix of the errors and its estimate, respectively. Consequently we obtain the GLS estimator of $\mbf{\b}$
\benn
\wh{\mbf{\b}}_{GLS}=(\textbf{X}'\wh{\mbf{\Omega}}^{-1}\textbf{X})(\textbf{X}'\wh{\mbf{\Omega}}^{-1}\mbf{y}).
\eenn
In order to obtain the generalized m.d.\,\,estimator, we try two different $\textbf{Q}$'s: $\textbf{Q}_{s}$ and $\textbf{Q}_{c}$ where
\benn
\textbf{Q}_{s}^{2}=\wh{\mbf{\Omega}}^{-1},\quad \textbf{Q}_{c}'\textbf{Q}_{c}=\wh{\mbf{\Omega}}^{-1}.
\eenn
We refer to the generalized m.d.\,\,estimators corresponding to $\textbf{Q}_{s}$ and $\textbf{Q}_{c}$
as GMD1 and GMD2 estimators, respectively.

In order to generate strongly mixing process for the dependent errors, the several restrictive conditions are required so that the mixing number $\al$ decays fast enough --- i.e., the assumption (a.8) is met. Withers (1981) proposed the upperbound and the decay rate of the mixing number $\al$. For the shake of completeness, we reproduce Theorem and Corollary 1 here.
\begin{lem}\lel{lem:Withers}
\textit{Let $\{\xi_{i}\}$ be independent r.v.s on R with characteristic functions $\{\phi_{i}\}$ such that}
\benn
(2\pi)^{-1}\max_{i}\int|\phi_{i}(t)|dt<\iny
\eenn
\textit{and}
\benn
\max_{i} E|\xi_{i}|^{\del}<\iny\quad\quad \textrm{for some }\del>0.
\eenn
\textit{Let $\{g_{v}:\, v=0,1,2,...\}$ be a sequence of complex numbers such that}
\benn
\Big\{S_{t}(\min(1,\del))\Big\}^{\max(1,\del)}\rightarrow 0 \quad\quad \textrm{as }t\rightarrow\iny
\eenn
\textit{where}
\benn
S_{t}(\lambda)=\sum_{v=t}^{\iny}|g_{v}|^{\lambda}.
\eenn
\textit{Assume that}
\benn
g_{v}=O(v^{-\kappa})\quad\textrm{ where }\kappa>1+\del^{-1}+\max(1,\del^{-1}).
\eenn
\textit{Then the sequence $\{ \vep_{n}:\,\vep_{n}= \sum_{v=0}^{\iny}g_{v}\xi_{n-v} \}$ is strongly mixing with mixing number $\al_{l}(k)=O(k^{-\eta})$ where}
\benn
\eta = (\kappa\del -\max(\del,1))(1+\del)^{-1}-1 > 0.
\eenn
\end{lem}
\noindent
To generate strongly mixing process by Lemma \r{lem:Withers}, we consider four independent $\xi_{i}$'s: normal, Laplace, logistic, and mixture of the two normals (MTN). Note that all the $\xi_{i}$'s have the finite second moments, and hence, we set $\del$ at 2. It can be easily seen that for any $\kappa>7$ we have $\eta>3$, and hence the assumption (a.8) is satisfied. Then for $\epsilon>0$
\ben\lel{eq:Gen_error}
\vep_{n}=\sum_{v=0}^{\iny}v^{-(7+\epsilon)}\xi_{n-v}
\een
satisfies the strongly mixing condition with $\al(k)=O(k^{-(3+\epsilon)})$. We let $\epsilon=0.5$, or equivalently, $\kappa=7.5$.

The $\xi$ has a Laplace distribution if its density function is
\benn
f_{La}(x):=(2s_{1})^{-1}\exp(-|x-\mu_{1}|/s_{1})
\eenn
while the density function of Logistic innovation is given by
\benn
f_{Lo}(x):= s_{2}^{-1}\exp(-|x-\mu_{2}|/s_{2})/(1+\exp(-|x-\mu_{2}|/s_{2}))^{2}.
\eenn
When we generate $\{\xi_{i}\}_{i=1}^{n}$, we set mean of normal, Laplace, and logistic innovations at 0 (i.e., $\mu_{1}=\mu_{2}=0$) since we assumed the $\vep$, the sum of $\xi_{i}$'s, is symmetric. We set the standard deviation of normal $\xi$ at 2 while both $s_{1}$ and $s_{2}$ are set at 5 for Laplace and logistic, respectively. For MTN, we consider $(1-\epsilon)N(0,2^2)+\epsilon N(0,10^2)$ where $\epsilon = 0.1$. In each $\xi$, we subsequently generate $\{\vep_{i}\}_{i=1}^{n}$ using (\r{eq:Gen_error}).

Next, we set the true $\mbf{\b}=(-2, 3, 1.5, -4.3)'$, i.e., $p=4$. For each $k=2,3,4$, we obtain $\{x_{ik}\}_{i=1}^{n}$ in \eqref{eq:Model1} as a random sample from the uniform distribution on $[0, 50]$; $\{y_{i}\}_{i=1}^{n}$ is subsequently generated using models (\r{eq:Model1}). We estimate $\mbf{\b}$ by the generalized m.d.\,\,and the GLS methods. We report empirical bias, standard error (SE), and mean squared error (MSE) of these estimators. We use the Lebesgue integrating measure, i.e., $H(y)\equiv y$. To obtain the generalized m.d.\,\,estimators, the author used \texttt{R} package \texttt{KoulMde}. The package is available from Comprehensive \texttt{R} Archive Network (CRAN) at \url{https://cran.r-project.org/web/packages/KoulMde/index.html}. Table 1 and 2 report biases, SE's and MSE's of estimators for the sample sizes 50 and 100, each repeated 1,000 times. The author used High Performance Computing Center (HPCC) to accelerate the simulations. All of the simulations were done in the \texttt{R}-3.2.2.

\begin{table}[htbp!]
\begin{center}
\begin{tabular}{c  c c c  c c c  c c c  c }
\hline
 &  & \multicolumn{3}{c}{GLS} & \multicolumn{3}{c}{GMD1} & \multicolumn{3}{c}{GMD2} \\
\cline{3-11}
 &  &  bias & SE & MSE & bias  & SE & MSE & bias & SE & MSE\\
\hline
\multirow{3}{*}{N}& $\b_{1}$&0.0011  &0.0798 &0.0064 &0.0014  &0.0797 &0.0064 &8e-04   &0.0805 &0.0065 \\
& $\b_{2}$&0.0036  &0.0753 &0.0057 &0.0032  &0.0756 &0.0057 &0.0036  &0.0761 &0.0058 \\
& $\b_{3}$&-0.0023 &0.0784 &0.0062 &-0.0024 &0.0784 &0.0062 &-0.0023 &0.0787 &0.0062 \\
& $\b_{4}$&0.0014  &0.0772 &0.006  &0.0016  &0.0773 &0.006  &0.0013  &0.0775 &0.006  \\
\hline
\multirow{3}{*}{La} &$\b_{1}$&-3e-04  &0.1085 &0.0118 &-3e-04  &0.107  &0.0115 &-6e-04  &0.1073 &0.0115\\
&$\b_{2}$&-0.0011 &0.1145 &0.0131 &-0.0011 &0.1137 &0.0129 &-0.001  &0.1147 &0.0132\\
&$\b_{3}$&-0.0011 &0.1129 &0.0127 &-0.001  &0.1121 &0.0126 &-0.0011 &0.1127 &0.0127 \\
&$\b_{4}$&7e-04   &0.1193 &0.0142 &7e-04   &0.119  &0.0142 &4e-04   &0.1194 &0.0143 \\
\hline
\multirow{3}{*}{Lo} &$\b_{1}$&-0.0111 &0.1438 &0.0208 &-0.0113 &0.1429 &0.0205 &-0.0108 &0.144  &0.0209\\
&$\b_{2}$&-0.0034 &0.1516 &0.023  &-0.0033 &0.1513 &0.0229 &-0.0033 &0.1515 &0.023\\
&$\b_{3}$&-0.0027 &0.1465 &0.0215 &-0.002  &0.1461 &0.0213 &-0.0024 &0.1465 &0.0215 \\
&$\b_{4}$&0.003   &0.1485 &0.0221 &0.0027  &0.1481 &0.0219 &0.0029  &0.1478 &0.0218 \\
\hline
\multirow{3}{*}{M}  &$\b_{1}$&-0.0024 &0.1005 &0.0101 &-0.0023 &0.0993 &0.0099 &-0.0027 &0.0996 &0.0099\\
&$\b_{2}$&0.0059  &0.1076 &0.0116 &0.0054  &0.1063 &0.0113 &0.0057  &0.1069 &0.0115\\
&$\b_{3}$&-0.002  &0.1035 &0.0107 &-0.0016 &0.1024 &0.0105 &-0.002  &0.1027 &0.0105 \\
&$\b_{4}$&-0.001  &0.1105 &0.0122 &-0.0013 &0.1098 &0.0121 &-0.001  &0.1097 &0.012  \\
\hline
\end{tabular}
\end{center}
\begin{tablenotes}
      \tiny
      \item $\dagger$ N, La, Lo, and M denote normal, Laplace, logistic and MTN, respectively.
\end{tablenotes}
\caption{Bias, SE, and MSE of estimators with $n=50$.}\lel{table:Panel-N105}
\end{table}

\begin{table}[htbp!]
\begin{center}
\begin{tabular}{c  c c c  c c c  c c c  c }
\hline
 &  & \multicolumn{3}{c}{GLS} & \multicolumn{3}{c}{GMD1} & \multicolumn{3}{c}{GMD2} \\
\cline{3-11}
 &  &  bias & SE & MSE & bias  & SE & MSE & bias & SE & MSE\\
\hline
\multirow{3}{*}{N}& $\b_{1}$&-3e-04  &0.0518 &0.0027 &0       &0.0522 &0.0027 &-5e-04  &0.0522 &0.0027\\
& $\b_{2}$&0.0014  &0.049  &0.0024 &0.0013  &0.0491 &0.0024 &0.0013  &0.0494 &0.0024\\
& $\b_{3}$&-0.0032 &0.0498 &0.0025 &-0.0031 &0.0499 &0.0025 &-0.0032 &0.0502 &0.0025\\
& $\b_{4}$&2e-04   &0.0495 &0.0024 &2e-04   &0.0496 &0.0025 &2e-04   &0.0496 &0.0025\\
\hline
\multirow{3}{*}{La} &$\b_{1}$&0.0013  &0.0731 &0.0053 &0.0018  &0.0725 &0.0053 &0.0015  &0.0725 &0.0053\\
&$\b_{2}$&-0.0027 &0.0703 &0.0049 &-0.0026 &0.0695 &0.0048 &-0.0025 &0.0697 &0.0049\\
&$\b_{3}$&-0.0039 &0.0715 &0.0051 &-0.004  &0.0712 &0.0051 &-0.0037 &0.071  &0.0051\\
&$\b_{4}$&-0.0015 &0.0672 &0.0045 &-0.0016 &0.0666 &0.0044 &-0.0018 &0.0668 &0.0045\\
\hline
\multirow{3}{*}{Lo} &$\b_{1}$&-5e-04 &0.0913 &0.0083 &-7e-04  &0.0907 &0.0082 &-2e-04  &0.0914 &0.0083\\
&$\b_{2}$&-0.001 &0.0925 &0.0086 &-0.0011 &0.0921 &0.0085 &-0.0011 &0.0926 &0.0086\\
&$\b_{3}$&0.0072 &0.0932 &0.0087 &0.0073  &0.0933 &0.0088 &0.0072  &0.0933 &0.0088\\
&$\b_{4}$&-4e-04 &0.0928 &0.0086 &-2e-04  &0.0929 &0.0086 &-3e-04  &0.093  &0.0087\\
\hline
\multirow{3}{*}{M}  &$\b_{1}$&-0.0029 &0.0684 &0.0047 &-0.0032 &0.067  &0.0045 &-0.0031 &0.0674 &0.0046\\
&$\b_{2}$&-8e-04  &0.069  &0.0048 &-0.0012 &0.0682 &0.0046 &-9e-04  &0.0679 &0.0046\\
&$\b_{3}$&5e-04   &0.0707 &0.005  &4e-04   &0.07   &0.0049 &3e-04   &0.0698 &0.0049\\
&$\b_{4}$&0.001   &0.0676 &0.0046 &0.0011  &0.0667 &0.0044 &0.0013  &0.0671 &0.0045\\
\hline
\end{tabular}
\end{center}
\caption{Bias, SE, and MSE of estimators with $n=100$.}\lel{table:Panel-N105}
\end{table}

As we expected, both biases and SE's of all estimators decrease as $n$ increases. First, we consider the normal $\xi$'s. When $\xi$'s are normal, the GLS and GMD1 estimators display the best performance; GLS and GMD1 show similar biases, SE's, and hence MSE's. GMD2 estimators show slightly worse performance than aforementioned ones; they display similar or smaller bias --- e.g. estimators corresponding to $n=50$ and $\b_{1}, \b_{3},\b_{4}$ --- while they always have larger SE's which in turn cause larger MSE's. Therefore, we conclude that GLS and GMD1 show similar performance to each other but better one than GMD2   when $\xi$'s are normal.

For non-Gaussian $\xi$'s, we come up with a different conclusion: the GMD1 estimators outperform all other estimators while The GLS and GMD2 estimators display the similar performance. Note that weighing the merits of the GLS, the GMD1, and the GMD2 estimators in terms of bias is hard. For example, for the Laplace $\xi$ when $n=50$, the GLS and GMD1 estimators of all $\b_{i}$'s show the almost same biases; the GMD2 estimator of $\b_{1}$ ($\b_{4}$) show smaller (larger) bias than the GLS and the GMD1 estimators. When we consider the SE, the GMD1 estimators display the least SE's regardless of $n$'s and $\xi$'s. The GLS and the GMD2 estimators show somewhat similar SE's when $\xi$ is Laplace or logistic; however, the GMD2 estimators have smaller SE's than the GLS ones when $\xi$ is MTN. As a result, the GMD1 estimators display the least MSE for all non-Gaussian $\xi$'s and $n$'s; the GMD2 and the GLS --- corresponding to Laplace or logistic $\xi$'s --- show similar MSE's while the GMD2 estimators show smaller MSE than the GLS ones when $\xi$ is MTN.

\appendix
\setcounter{section}{1}
\section*{Appendix}
\textbf{Proof of Theorem \r{thm:1}.} Section 5.5 of Koul (2002) illustrates (\r{eq:thm1}) holds for independent errors. Proof of the theorem, therefore, will be similar to the one of Theorem 5.5.1 in that section. Define for $k=1,2,...,p$, $\mbf{u}\in \R^{p}, y\in \R$,
\benr\lel{eq:defineJk}
J_{k}(y,\mbf{u})&:=& \sti d_{ik}\,F_{i}(y + \mbf{q}_{i}'\textbf{X}\textbf{A}\mbf{u}),\quad
Y_{k}(y,\mbf{u}) :=\sti d_{ik}\,I\big(\mbf{q}_{i}'\mbf{\vep}\leq y + \mbf{q}_{i}'\textbf{X}\textbf{A}\mbf{u}  \big),\\
W_{k}(y,\mbf{u}) &:=& Y_{k}(y,\mbf{u}) - J_{k}(y,\mbf{u}).\nonumber
\eenr
Rewrite
\benr\lel{eq:Toft}
\,\,\,\,\,\,\,\,\,\cL(\mbf{\b}+ \textbf{A}\mbf{u};\,\textbf{Q}) &=& \sum_{k=1}^{p}\int\; \Big[ \left\{ W_{k}(y,\mbf{u}) - W_{k}(y,\mbf{0})\right\} + \left\{ W_{k}(-y,\mbf{u}) - W_{k}(-y,\mbf{0})\right\} \\
&&\quad +  \{ (J_{k}(y,\mbf{u}) - J_{k}(y,\mbf{0})) - \sti d_{ik}\mbf{q}_{i}'\textbf{X}\textbf{A}\mbf{u}f_{i}(y)  \}\nonumber\\
&&\quad + \{ (J_{k}(-y,\mbf{u}) - J_{k}(-y,\mbf{0})) - \sti d_{ik}\mbf{q}_{i}'\textbf{X}\textbf{A}\mbf{u}f_{i}(y) \} \nonumber \\
&&\quad +   \{ U_{k}(y,\mbf{\b})+2\sti d_{ik}\mbf{q}_{i}'\textbf{X}\textbf{A}\mbf{u}f_{i}(y) \}\Big]^{2}\,d H(y). \nonumber
\eenr
where $\mbf{0}=(0,0,...,0)'\in\R^{p}$. Note that the last term of the integrand is the $k$th coordinate of  $\mbf{U}(y,\mbf{\b};\,\textbf{Q}) + 2\mbf{\Sigma}_{\textbf{DA}}(y)\textbf{A}^{-1}(\mbf{b}-\mbf{\b})$ vector in $\cL^{*}(\mbf{b};\,\textbf{Q})$. If we can show that suprema of $L_{H}^{2}$ norms of the first four terms of the integrand are $o_{p}(1)$, then applying Cauchy-Schwarz (C-S) inequality on the cross product terms in \eqref{eq:Toft} will complete the proof. Therefore to prove theorem it suffices to show that for all $k=1,2,...,p$
\ben\lel{eq:thm11}
E\sup \int\,\big|W_{k}(\pm y, \mbf{u})-W_{k}(\pm y,\mbf{0})\big|^{2} \,d H(y)=o(1),
\een
\ben\lel{eq:thm12}
\sup\int\,\big| (J_{k}(\pm y,\mbf{u}) - J_{k}(\pm y,\mbf{0})) - \sti d_{ik}\mbf{q}_{i}'\textbf{X}\textbf{A}\mbf{u}f_{i}(y)\big|^{2}\, d H(y)=o(1),
\een
\ben\lel{eq:thm13}
E\sup \int\, \big| U_{k}(y,\mbf{\b})  + 2\sti d_{ik}\mbf{q}_{i}'\textbf{X}\textbf{A}\mbf{u}f_{i}(y) \big|^{2}\, d H(y) =O(1).
\een
where sup is taken over $\|\mbf{u}\|\leq b$. Here we consider the proof of the case $+y$ only. The similar facts will hold for the case $-y$.

Observe that (A.2) implies
\ben\lel{eq:A4}
E\int U_{k}(y,\mbf{\b})^{2}dH(y)\le 2n\max_{1\le i\le n}\|\mbf{d}_{i}\|^2 \max_{1\le i\le n} \int (1-F_{i})\,dH<\iny.
\een
Therefore, (\r{eq:thm13}) immediately follows from (A.2) and (A.7). The proof of (\r{eq:thm12}) does not involve the dependence of errors, and hence, it is the same as the proof of (5.5.11) of Koul (2002). Thus, we shall prove (\r{eq:thm11}), thereby completing the proof of theorem.

To begin with let $J_{ku}^{\pm}(\cdot)$, $Y_{ku}^{\pm}(\cdot)$, and $W_{ku}^{\pm}(\cdot)$ denote $J_{k}(\cdot, \mbf{u})$, $Y_{k}(\cdot, \mbf{u})$ and $W_{k}(\cdot, \mbf{u})$ in \eqref{eq:defineJk} when $d_{ik}$ is replaced with $d_{ik}^{\pm}$ so that $J_{k}=J_{k}^{+}-J_{k}^{-}$, $Y_{k}=Y_{k}^{+}-Y_{k}^{-}$, and $W_{k}=W_{k}^{+}-W_{k}^{-}$. Define for $x\in \R^{p},\,u\in \R^{p}$, $y\in \R$,
\benrr
p_{i}(y,\mbf{u}; \textbf{X})&:=& F_{i}(y+\mbf{q}_{i}'\textbf{X}\textbf{A}\mbf{u})-F(y),\\
B_{ni}&:=&  I\big( \mbf{q}_{i}'\mbf{\vep}\leq y + \mbf{q}_{i}'\textbf{X}\textbf{A}\mbf{u}\big) -I\big(\mbf{q}_{i}'\mbf{\vep}\leq y \big) - p_{i}(y,\mbf{u}; \textbf{X}).
\eenrr
Rewrite
\benrr
W_{ku}^{\pm}-W_{k0}^{\pm} = \sti d_{ik}^{\pm}\,\big\{ I\big( \mbf{q}_{i}'\mbf{\vep}\leq y + \mbf{q}_{i}'\textbf{X}\textbf{A}\mbf{u}\big) -I\big(\mbf{q}_{i}'\mbf{\vep}\leq y \big) - p_{i}(y,\mbf{u}; \textbf{X})  \big\}.
\eenrr
Note that
\ben\lel{eq:Bni}
E\,B_{ni}^{2} \leq  F_{i}(y+ \theta_{i}\|\mbf{u}\|)-F_{i}(y)
\een
Recall a lemma from Deo(1973).
\begin{lem}\lel{lemma:Deo}
Suppose for each $n\geq 1$, $\{\xi_{nj},\,1\leq j\leq n\}$ are strongly mixing random variables with mixing number $\alpha_{n}$. Suppose $X$ and $Y$ are two random variables respectively measurable with respect to $\sigma\{\xi_{n1},...,\xi_{nk}\}$ and $\sigma\{\xi_{nk+m},...,\xi_{nn}\}$, $1\leq m,\,m+k\leq n$. Assume $p,q$ and $r$ are such that $p^{-1}+q^{-1}+r^{-1}$, and $\|X\|_{p}\leq\iny$ and $\|Y\|_{q}\leq\iny$. Then for each $1\leq m,\,k+m\leq n$
\ben
|E(XY)-E(X)E(Y)|\leq 10\cdot\alpha_{n}^{1/r}(m)\|X\|_{p}\|Y\|_{q}.
\een
Consequently if $\|X\|_{\iny}=B<\iny$ then for $q>1$ and each $1\leq m,\,k+m\leq n$
\ben
|E(XY)-E(X)E(Y)|\leq 10\cdot\alpha_{n}^{1-1/q}(m)\|Y\|_{q}.
\een
\end{lem}
\noindent
In addition, consider following lemma.
\begin{lem}\lel{lemma:Infinite_Test}
For $1<r<3$,
\ben
n^{-1}\sum_{i=1}^{n-1}\sum_{k=1}^{n-i}\alpha^{1/r}(k)=O(1). \lel{eq:Infinite_Test}
\een
\end{lem}
\noindent
\textbf{Proof}. For given $r$, let $p$ such that $1/r+1/p=1$. Note that
\benn
\frac{p}{r} = \frac{1}{r-1}>\frac{1}{2}.
\eenn
Therefore, by H$\ddot{\textrm{o}}$lder's inequality with $p$ and $r$, we have
\benr\lel{eq:lemma_InfTest}
n^{-1}\sum_{i=1}^{n-1}\sum_{k=1}^{n-i}\alpha^{1/r}(k)
&\leq& \left(\sum_{k=1}^{n-1}\frac{(n-k)^{p}}{n^{p}}\cdot\frac{1}{k^{2p/r}} \right)^{1/p}  \left(\sum_{k=1}^{n-1}k^{2}\alpha(k)\right)^{1/r} <\iny.
\eenr
The last inequality follows from the assumption (A.8.6), thereby completing the proof of lemma. \qed

Now, we consider the cross product terms of $E \big|W_{ku}-W_{k0}\big|_{H}^{2}$.
\benr
&&\Big|E\, \int \sti\sum_{j=i+1}^{n} \left[\, d_{ik}d_{jk}\big\{ I\big( \mbf{q}_{i}'\mbf{\vep}\leq y + \mbf{q}_{i}'\textbf{X}\textbf{A}\mbf{u}\big) -I\big(\mbf{q}_{i}'\mbf{\vep}\leq y \big) - p_{i}(y,\mbf{u}; \textbf{X}) \big\} \right. \lel{eq:Covariance}\\
&& \quad\quad \left. \times \big\{ I\big( \mbf{q}_{j}'\mbf{\vep}\leq y + \mbf{q}_{j}'\textbf{X}\textbf{A}\mbf{u}\big) -I\big(\mbf{q}_{j}'\mbf{\vep}\leq y \big) - p_{j}(y,\mbf{u}; \textbf{X}) \big\} \right]\, d H(y)\Big|\nonumber\\
&\leq& \sti\sum_{j=i+1}^{n} d_{ik}d_{jk}\,\int \left| E\, B_{ni}\,B_{nj}   \right|\,d H\nonumber\\
&\leq& 10\sti\sum_{j=i+1}^{n} d_{ik}d_{jk}\,\cdot \alpha^{1/2}(j-i)\cdot \int  \|B_{nj}\|_{2} \,d H\nonumber\\
&\leq& 10b^{1/2}\{n\max_{i}d_{ik}^{2}\}\cdot\{\max_{i}\theta_{i}\}^{1/2}\cdot n^{-1}\sum_{i=1}^{n-1}\sum_{m=1}^{n-i} \alpha^{1/2}(m) \int  f_{1}^{1/2}  \,d H \ra 0.\nonumber
\eenr
The second inequality follows from Lemma \r{lemma:Deo}, and the convergence to zero follows from the Lemma \r{lemma:Infinite_Test} with $r=2$, (A.1), and (A.7). Consequently, by Fubini's Theorem together with (A.3), we obtain, for every fixed $\|\mbf{u}\|\leq b$,
\benr\lel{eq:DiniW}
\limsup_{n\ra\iny}E\,|W_{ku}-W_{k0}|^{2}_{H} &\leq & \limsup_{n\ra\iny} \int \sti d_{ik}^{2}\big|F_{i}(y+\mbf{q}_{i}'\textbf{X}\textbf{A}\mbf{u})-F(y) \big|  \,d H(y)\\
&\leq&\limsup_{n\ra\iny} \{n\max_{i}d_{ik}^{2}\} \int_{-a_{n}}^{a_{n}}\,\int f_{i}(y+s) d H(y)ds\nonumber\\
&=&0,\nonumber
\eenr
where $a_{n}=b\max_{i}\theta_{i}\ra0$.

To complete the proof of (\r{eq:thm11}), it suffices to show that for all $\epsilon>0$, there exists a $\delta>0$ such that for all $\mbf{v}\in \R^{p}$, $\|\mbf{u}-\mbf{v}\|\leq \delta$,
\ben\lel{eq:thm111}
\limsup_{n\ra\iny} E\sup_{\|u-v\|\leq \delta}|\cK_{ku}-\cK_{kv}|\leq \epsilon,
\een
\noindent
where
\benn
\cK_{ku}:=|W_{ku}-W_{k0}|_{H}^{2},\quad \mbf{u} \in\R^{p},\quad 1\leq k\leq p.
\eenn
(\r{eq:thm111}) follows from (5.5.5) of Koul (2002), thereby completing the proof of theorem.\qed

\noindent
\textbf{Proof of Corollary \r{cor:1}.} The proof of (\r{eq:cor1}) for independent errors can again be found in the section 5.5 of Koul (2002). The difference between the proof in the section 5.5 and one here arises only in the part which involves the dependence of the error. Thus, we present only the proof of an analogue of (5.5.27) in Koul (2002). Let
\benr
L_{k}&:=& \int\,\Big[ W_{k}(y,\mbf{0})+ W_{k}(-y,\mbf{0})+ \big\{ J_{k}(y,\mbf{0})+J_{k}(-y,\mbf{0}) - \sum_{i=1}^{n}d_{ik}\big\}\Big] f_{1}(y)\,d H(y)\nonumber\\
&=& \int\,\Big[\sti d_{ik}\left\{  I\big(\mbf{q}_{i}'\mbf{\vep}\leq y\big) - I\big(-\mbf{q}_{i}'\mbf{\vep} < y\big)   \right\}  \Big]f_{1}(y)\,d H(y),\nonumber\\
\mbf{L}&:=& (L_{1},...,L_{p}).\nonumber
\eenr
Note that $L_{k}=\int U_{k}(y,\mbf{\b};\,\textbf{Q})f_{1}(y)d H(y)$. By the symmetry of $H$ and Fubini's theorem, we obtain
\benn
E\int\left\{  I\big(\mbf{q}_{i}'\mbf{\vep}\leq y\big) - I\big(-\mbf{q}_{i}'\mbf{\vep} < y\big)   \right\}^{2}d H(y)=4\int_{0}^{\iny}(1-F_{i})\,d H,\quad 1\le i\le n.
\eenn
In addition, Lemma \r{lemma:Deo} yields, for $j>i$,
\benr
&&E\int\left\{  I\big(\mbf{q}_{i}'\mbf{\vep}\leq y\big) - I\big(-\mbf{q}_{i}'\mbf{\vep} < y\big)   \right\}\left\{  I\big(\mbf{q}_{j}'\mbf{\vep}\leq y\big) - I\big(-\mbf{q}_{j}'\mbf{\vep} < y\big)   \right\}d H(y)\nonumber\\
&\leq & 20\sqrt{2} \,\,\alpha^{1/2}(j-i)\,\,\max_{1\le i\le n}\int_{0}^{\iny}(1-F_{i})^{1/2}d H.\nonumber
\eenr
Together with the fact that $(1-F_{i})\leq (1-F_{i})^{1/2}$, by (A.1), (A.2), (A.7), and Lemma \r{lemma:Infinite_Test}, we obtain, for some $0<C<\iny$,
\benr
E\|\mbf{L}\|^{2}&\leq& \sum_{k=1}^{p}|f_{1}|_{H}^{2}\Big\{ 4\sum_{i=1}^{n}d_{ik}^{2}\int_{0}^{\iny}(1-F_{i})\,d H + 40\sqrt{2}\,\,n\max_{1\leq i\leq n}d_{ik}^{2} \nonumber\\
& & \quad\quad\quad\quad\quad\quad\quad  \times n^{-1}\sum_{i=1}^{n}\sum_{j=i+1}^{n}\alpha^{1/2}(j-i)\,\,\max_{1\leq i\leq n}\int_{0}^{\iny}(1-F_{i})^{1/2}\,d H \Big\}\nonumber\\
&<&C\,p\,|f_{1}|_{H}^{2}\,\max_{1\leq i\leq n}\int_{0}^{\iny}(1-F_{i})^{1/2}\,d H. \nonumber
\eenr
Using $E\,L_{k}=0$ for $k=1,...,p+1$ and Chebyshev inequality,
for all $\epsilon>0$ there exists $N_{1}$ and $c_{\epsilon}$ such that
\ben\lel{eq:Cor13}
P\big(\|\mbf{L}\|\leq c_{\epsilon}\big) \geq  1-\frac{C\,p\,|f_{1}|_{H}^{2}\,\max_{1\leq i\leq n}\int_{0}^{\iny}(1-F_{i})^{1/2}\,d H  }{c_{\epsilon}} \geq 1- \epsilon/2,\quad n\geq N_{1}.
\een
The rest of the proof will be the same as the proof of Lemma 5.5.4 of Koul (2002).

\edt

\section*{Appendix}
\noindent
\textbf{Proofs of Section \r{asnb}}\\
\noindent
\textbf{Proof of Theorem \r{thm:1}.} Section 5.5 of Koul (2002) illustrates the claim holds for independent errors. Proof of the theorem, therefore, will be similar to the one of Theorem 5.5.1 in that section. Define for $k=1,2,...,p+1$, $u\in \R^{p+1}, y\in \R$,
\benr\lel{eq:defineJk}
J_{k}(y,u)&:=& \sti  d_{ik}\,F(y + u'Ax_{i}),\quad
Y_{k}(y,u) :=\sti  d_{ik}\,I\big(\vep_{i}\leq y + u'Ax_{i}  \big),\\
W_{k}(y,u) &:=& Y_{k}(y,u) - J_{k}(y,u),\nonumber
\eenr
where $x_{i}\in \R^{p+1},\, i=1,2,...,n$ are in the model \eqref{eq:Model1}, and $d_{ik}=a_{k}'x_{i}$. Note that $A\sti  d_{ik}x_{i}=e_{k}$, where $e_{k}\in\R^{p+1}$ is elementary vector whose $k$th entry is 1. Therefore, rewrite
\benr\lel{eq:Toft}
\,\,\,\,\,\,\,\,\,T(\b+ Au) &=& \sum_{k=1}^{p+1}\int\; \Big[ \left\{ W_{k}(y,u) - W_{k}(y,0)\right\} + \left\{ W_{k}(-y,u) - W_{k}(-y,0)\right\} \\
&&\quad +  \{ (J_{k}(y,u) - J_{k}(y,0)) - u'e_{k}\,f(y)  \}\nonumber\\
&&\quad + \{ (J_{k}(-y,u) - J_{k}(-y,0)) - u'e_{k}\,f(y) \} \nonumber \\
&&\quad +   \{ U_{k}(y,\b)+2u'e_{k}\,f(y) \}\Big]^{2}\,d H(y). \nonumber
\eenr
Note that the last term of the integrand is the $k$th coordinate of  $U(y,\b) + 2A^{-1}(t-\b)f(y)$ vector in $Q(t)$. If we can show that suprema of $L_{H}^{2}$ norms of the first four terms of the integrand are $o_{p}(1)$, then applying Cauchy-Schwarz (C-S) inequality on the cross product terms in \eqref{eq:Toft} will complete the proof. Therefore to prove theorem it suffices to show that for all $k=1,2,...,p+1$,
\benr
&&E\sup \int\,\big|W_{k}(\pm y,u)-W_{k}(\pm y,0)\big|^{2} \,d H(y)=o(1),\lel{eq:thm11}\\
&&\sup\int\,\big| (J_{k}(\pm y,u) - J_{k}(\pm y,0)) - u'e_{k}\,f(y)\big|^{2}\, d H(y)=o(1),\lel{eq:thm12}\\
&&E\sup \int\, \big| U_{k}(y,\b)  + 2u'e_{k}\,f(y)\big|^{2}\, d H(y) =O(1)\lel{eq:thm13}.
\eenr
where sup is taken over $\|u\|\leq b$.

First consider (\r{eq:thm13}). By the symmetry of $F$ and $H$ and Fubini's theorem, we obtain
\benn
E\int\left\{  I\big(\vep_{i}\leq y\big) - I\big(-\vep_{i} < y\big)   \right\}^{2}d H(y)=4\int_{0}^{\iny}(1-F)\,d H.
\eenn
In addition, Lemma \r{lemma:Deo} yields, for $j>i$,
\benr
&&E\int\left\{  I\big(\vep_{i}\leq y\big) - I\big(-\vep_{i} < y\big)   \right\}\left\{  I\big(\vep_{j}\leq y\big) - I\big(-\vep_{j} < y\big)   \right\}d H(y)\nonumber\\
&\leq & 20\sqrt{2} \,\,\alpha_{l}^{1/2}(j-i)\int_{0}^{\iny}(1-F)^{1/2}d H.\nonumber
\eenr
Together with the fact that $(1-F)\leq (1-F)^{1/2}$, by (A.1), (A.2), (A.7), and Lemma \r{lemma:Infinite_Test}, we obtain
\benr\lel{eq:EUL}
E\int U_{k}(y,\b)^{2}dH(y) &\leq&  4 n\max_{1\leq i\leq n}d_{ik}^{2}\int_{0}^{\iny}(1-F)\,d H + 40\sqrt{2}\,\,n\max_{1\leq i\leq n}d_{ik}^{2} \\
& & \quad\quad\quad\quad\quad  \times n^{-1}\sum_{i=1}^{n}\sum_{j=i+1}^{n}\alpha_{l}^{1/2}(j-i)\int_{0}^{\iny}(1-F)^{1/2}\,d H \nonumber\\
&<&\iny. \nonumber
\eenr
Therefore, (\r{eq:thm13}) immediately follows.

The proof of (\r{eq:thm12}) does not involve the dependence of errors, and hence, it is the same as the proof of (5.5.11) of Koul (2002). Thus, we shall prove (\r{eq:thm11}). Here we consider the proof for the case $+y$ only. The similar facts will hold for the case $-y$.

To begin with let $J_{ku}^{\pm}(\cdot)$, $Y_{ku}^{\pm}(\cdot)$, and $W_{ku}^{\pm}(\cdot)$ denote $J_{k}(\cdot, u)$, $Y_{k}(\cdot, u)$ and $W_{k}(\cdot, u)$ in \eqref{eq:defineJk} when $d_{ik}$ is replaced with $d_{ik}^{\pm}$ so that $J_{k}=J_{k}^{+}-J_{k}^{-}$, $Y_{k}=Y_{k}^{+}-Y_{k}^{-}$, and $W_{k}=W_{k}^{+}-W_{k}^{-}$. Define for $x\in \R^{p+1},\,u\in \R^{p+1}$, $y\in \R$,
\benrr
p(y,u; x)&:=& F(y+u'Ax)-F(y),\\
B_{ni}&:=&  I\big(\vep_{i}\leq y + u'Ax_{i}\big) -I\big(\vep_{i}\leq y \big) - p(y,u; x_{i}).
\eenrr
Rewrite
\benrr
W_{ku}^{\pm}-W_{k0}^{\pm} = \sti   d_{ik}^{\pm}\,\big\{I\big(\vep_{i}\leq y + u'Ax_{i}\big)  -I\big(\vep_{i}\leq y \big) - p(y,u; x_{i})     \big\}.
\eenrr
Note that
\ben\lel{eq:Bni}
E\,B_{ni}^{2} \leq  F(y+ \theta_{i}\|u\|)-F(y)
\een
Recall a lemma from Deo(1973).
\begin{lem}\lel{lemma:Deo}
Suppose for each $n\geq 1$, $\{\xi_{nj},\,1\leq j\leq n\}$ are strongly mixing random variables with mixing number $\alpha_{l}$. Suppose $X$ and $Y$ are two random variables respectively measurable with respect to $\sigma\{\xi_{n1},...,\xi_{nk}\}$ and $\sigma\{\xi_{nk+m},...,\xi_{nn}\}$, $1\leq m,\,m+k\leq n$. Assume $p,q$ and $r$ are such that $p^{-1}+q^{-1}+r^{-1}$, and $\|X\|_{p}\leq\iny$ and $\|Y\|_{q}\leq\iny$. Then for each $1\leq m,\,k+m\leq n$
\benn
|E(XY)-E(X)E(Y)|\leq 10\cdot\alpha_{l}^{1/r}(m)\|X\|_{p}\|Y\|_{q}.
\eenn
Consequently if $\|X\|_{\iny}=B<\iny$ then for $q>1$ and each $1\leq m,\,k+m\leq n$
\benn
|E(XY)-E(X)E(Y)|\leq 10\cdot\alpha_{l}^{1-1/q}(m)\|Y\|_{q}.
\eenn
\end{lem}
\noindent
In addition, consider following lemma.
\begin{lem}\lel{lemma:Infinite_Test}
For $1<r<3$,
\benn
n^{-1}\sum_{i=1}^{n-1}\sum_{k=1}^{n-i}\alpha_{l}^{1/r}(k)=O(1).
\eenn
\end{lem}
\noindent
\textbf{Proof}. Fix an $1<r<3$. Let $p$ such that $1/r+1/p=1$. Note that
\benn
\frac{p}{r} = \frac{1}{r-1}>\frac{1}{2}.
\eenn
Therefore, by H$\ddot{\textrm{o}}$lder's inequality, we have
\benr\lel{eq:lemma_InfTest}
n^{-1}\sum_{i=1}^{n-1}\sum_{k=1}^{n-i}\alpha_{l}^{1/r}(k)
&\leq& \left(\sum_{k=1}^{n-1}\frac{(n-k)^{p}}{n^{p}}\cdot\frac{1}{k^{2p/r}} \right)^{1/p}  \left(\sum_{k=1}^{n-1}k^{2}\alpha_{l}(k)\right)^{1/r} <\iny.
\eenr
The last inequality follows from the assumption (A.8.6), thereby completing the proof of the lemma.

Now, we consider the cross product terms of $E \big|W_{ku}-W_{k0}\big|_{H}^{2}$. Then
\benr
&&\Big|E\, \int \sti  \sum_{j=i+1}^{n} \left[\, d_{ik}d_{jk}\big\{ I\big(\vep_{i}\leq y + u'Ax_{i}\big) -I\big(\vep_{i}\leq y \big) - p(y,u; x_{i}) \big\} \right. \nonumber\\
&& \quad\quad \left. \times \big\{ I\big(\vep_{j}\leq y + u'Ax_{j}\big) -I\big(\vep_{j}\leq y \big) - p(y,u; x_{j}) \big\} \right]\, d H(y)\Big|\nonumber\\
&\leq& \sti  \sum_{j=i+1}^{n} d_{ik}d_{jk}\,\int \left| E\, B_{ni}\,B_{nj}   \right|\,d H\nonumber\\
&\leq& 10b^{1/2}\{n\max_{i}d_{ik}^{2}\}\cdot\{\max_{i}\theta_{i}\}^{1/2}\cdot n^{-1}\sum_{i=1}^{n-1}\sum_{m=1}^{n-i} \alpha_{l}^{1/2}(m) \int  f^{1/2}  \,d H \ra 0.\nonumber
\eenr
The second inequality follows from Lemma \r{lemma:Deo}, and the convergence to zero follows from the Lemma \r{lemma:Infinite_Test} with $r=2$, (A.1), and (A.7). Consequently, by Fubini's Theorem together with (A.3), we obtain, for every fixed $\|u\|\leq b$,
\benr\lel{eq:DiniW}
\limsup_{n\ra\iny}E\,|W_{ku}-W_{k0}|^{2}_{H} &\leq & \limsup_{n\ra\iny}  \int_{-a_{n}}^{a_{n}}\,\int f(y+s) d H(y)ds =0,\nonumber
\eenr
where $a_{n}=b\max_{i}\theta_{i}\ra0$.

To complete the proof of (\r{eq:thm11}), it suffices to show that for all $\epsilon>0$, there exists a $\delta>0$ such that for all $v\in \R^{p+1}$, $\|u-v\|\leq \delta$,
\ben\lel{eq:thm111}
\limsup_{n\ra\iny} E\sup_{\|u-v\|\leq \delta}|\cK_{ku}-\cK_{kv}|\leq \epsilon,
\een
\noindent
where
\benn
\cK_{ku}:=|W_{ku}-W_{k0}|_{H}^{2},\quad u\in\R^{p+1},\quad 1\leq k\leq p+1.
\eenn
(\r{eq:thm111}) follows from (5.5.5) of Koul (2002), thereby completing the proof of theorem.

\noindent
\textbf{Proof of Corollary \r{cor:1}.} The proof of the claim for independent errors can again be found in the section 5.5 of Koul (2002). The difference between the proof in the section 5.5 and one here arises only in the part which involves the dependence of the error. Thus, we present only the proof of an analogue of (5.5.27) in Koul (2002). Let
\benr
L_{k}&:=& \int\,\Big[\sti  d_{ik}\left\{  I\big(\vep_{i}\leq y\big) - I\big(-\vep_{i} < y\big)   \right\}  \Big]f(y)\,d H(y),\nonumber\\
L&:=& (L_{1},...,L_{p+1}).\nonumber
\eenr
Note that $L_{k}=\int U_{k}(y,\b)f(y)d H(y)$. Therefore, using C-S inequality together with (\r{eq:EUL}), we obtain, for some $0<C<\iny$
\benn
E\|L\|^{2} \leq C(p+1)|f|_{H}^{2}\int_{0}^{\iny}(1-F)^{1/2}\,d H. \nonumber
\eenn
Using $E\,L_{k}=0$ for $k=1,...,p+1$ and Chebyshev inequality,
for all $\epsilon>0$ there exists $N_{1}$ and $c_{\epsilon}$ such that
\benn
P\big(\|L\|\leq c_{\epsilon}\big) \geq  1-\frac{C(p+1) |f|_{H}^{2}\,\int(1-F)^{1/2}d H  }{c_{\epsilon}} \geq 1- \epsilon/2,\quad n\geq N_{1}.
\eenn
The rest of the proof is the same as that of Lemma 5.5.4 of Koul (2002).

\noindent\\
\textbf{Proofs of Section \r{asnro}}\\
\noindent
\textbf{Proof of Theorem \r{thm:sec_rho_thm1}}. Note that
\benrr
\widehat{S}_{k}(y,r)&=&S_{k}(y,r)+\, n^{-1/2} \sti   \vep_{i-k} \Bigg[\Big\{\,I\big( \xi_{i} \leq y+Z_{i}'v/\sqrt{n} +\z_{ni}(u,v)\big)-I\big( \xi_{i} \leq y+Z_{i}'v/\sqrt{n} \big)\,\Big\} \\
&&\quad- \Big\{\,I\big( -\xi_{i} \leq y-Z_{i}'v/\sqrt{n} -\z_{ni}(u,v)\big)-I\big( -\xi_{i} < y-Z_{i}'v/\sqrt{n} \big)\,\Big\}\Bigg]\\
&& -\, n^{-1/2} \sti   \eta_{n,i-k}(u) \Big\{\,I\big(\xi_{i}  \leq y + Z_{i}'v/\sqrt{n} \big) - I\big(-\xi_{i} < y - Z_{i}'v/\sqrt{n} \big)\,\Big\} \\
&& -\, n^{-1/2} \sti   \eta_{n,i-k}(u) \Bigg[\Big\{\,I\big( \xi_{i} \leq y+Z_{i}'v/\sqrt{n} +\z_{ni}(u,v)\big)-I\big( \xi_{i} \leq y+Z_{i}'v/\sqrt{n} \big)\,\Big\}\\
&&\quad - \Big\{\,I\big( -\xi_{i} \leq y-Z_{i}'v/\sqrt{n} -\z_{ni}(u,v)\big)-I\big( -\xi_{i} < y-Z_{i}'v/\sqrt{n} \big)\,\Big\}\Bigg]\\
&=& S_{r} + R1_{uv}^{k}- R2_{uv}^{k} - R3_{uv}^{k},\quad (say).
\eenrr
In what follows, the index $k$ runs from 1 to $p+1$. In order to conserve the space, we write
$\eta_{u}^{ik}$, $\z_{iuv}$, etc. for $\eta_{n,i-k}(u)$, $\z_{ni}(u,v)$, etc. Similar to the way as that used in Theorem \r{thm:1}, in order to prove Theorem \r{thm:sec_rho_thm1}, it suffices to show that for $j=1,2,3$
\benn
\sup_{\|u\|\leq b,\|v\|\leq b} \big|Rj_{uv}^{k}\big|_{H}^{2}=o_{p}(1).
\eenn
\noindent
Consider $R1_{uv}^{k}$. First, define for $1\leq k\leq q$
\benrr
w_{k}(y,u,v)&:=& n^{-1}\sti  \vep_{i-k}I\big( \xi_{i} \leq y+Z_{i}'v/\sqrt{n} +\z_{ni}(u,v)\big),\\
\nu_{k}(y,u,v)&:=&n^{-1}\sti  \vep_{i-k}F_{\tiny{\xi}}(y+Z_{i}'v/\sqrt{n} +\z_{ni}(u,v)),\\
\cW_{k}(y,u,v)&:=& \sqrt{n}\big[w_{k}(y,u,v)-\nu_{k}(y,u,v)\big].
\eenrr
Also consider the following processes:
\benrr
\cT_{k}(y;u,v, \pm\delta)&:=&n^{-1/2}\sti  \vep_{i-k}I\big( \xi_{i} \leq y+Z_{i}'v/\sqrt{n} +\z_{ni}(u,v)\pm d_{ni}(\delta, u,v)\big),\\
m_{k}(y;u,v, \pm\delta)&:=& n^{-1/2}\sti  \vep_{i-k}F_{\tiny{\xi}}(y+Z_{i}'v/\sqrt{n} +\z_{ni}(u,v) \pm d_{ni}(\delta, u,v)),\\
\cZ_{k}(y;u,v, \pm\delta) &:=& \cT_{k}(y;u,v, \pm\delta) - m_{k}(y;u,v, \pm\delta).
\eenrr
Note that $d_{ni}(-\delta, u,v)\neq -d_{ni}(\delta, u,v)$. However, for the purpose of conserving the space, let
\benn
\cT_{k}(y;u,v, -\delta):=n^{-1/2}\sti  \vep_{i-k}I\big( \xi_{i} \leq y+Z_{i}'v/\sqrt{n} +\z_{ni}(u,v)-d_{ni}(\delta, u,v)\big).
\eenn
Similarly, we define $m_{k}(y;u,v, -\delta)$ and $\cZ_{k}(y;u,v, -\delta)$ by subtracting term $d_{ni}(\delta, u,v)$. Let $\vep_{i-k}^{+}=0\vee \vep_{i-k}$ and $\vep_{i-k}^{-}=\vep_{i-k}^{+}-\vep_{i-k}$. Let any process with $\pm$ superscript denote the one with $\vep_{i-k}$ replaced by $\vep_{i-k}^{\pm}$. Also drop the subscript $k$, and let $w_{uv},\,\nu_{uv}, \,\cW_{uv}$, etc. stand for $w_{k}(
\cdot,u,v),\,\nu_{k}(\cdot,u,v),\,\cW_{k}(\cdot,u,v)$, etc.

Now, consider $s\in \R^{p+1},t\in \R^{q}$ where $\|s-u \|\leq \delta$ and $\|t-v \|\leq \delta$. Then we obtain
\ben\lel{eq:strange}
Z_{i}'v/\sqrt{n} +\z_{ni}(u,v)-d_{ni}(\delta,u,v )\leq Z_{i}'t/\sqrt{n} +\z_{ni}(s,t)\leq Z_{i}'v/\sqrt{n} +\z_{ni}(u,v)+d_{ni}(\delta,u,v ).
\een
Note that $\cT(y;u,v, 0)=\sqrt{n}\,w_{uv}$ and $m(y;u,v, 0)=\sqrt{n}\,\nu_{uv}$. Monotonicity of indicator function together with (\r{eq:strange}) in turn implies
\benn
\cT^{\pm}(y;u,v, -\delta)-\cT^{\pm}(y;u,v, 0)\leq \sqrt{n}\,[w_{st}^{\pm}-w_{uv}^{\pm}]\leq \cT^{\pm}(y;u,v, \delta)-\cT^{\pm}(y;u,v, 0)
\eenn
for all $y\in \R$, $\|s\|\leq b$, $\|t\|\leq b$, $\|s-u \|\leq \delta$ and $\|t-v \|\leq \delta$. Using the fact that $a_{1}\leq a_{2}\leq a_{3}$ implies $|a_{2}|\leq |a_{1}|+|a_{3}|$, appropriate centering of $\cT^{\pm}$ yields
\benr\lel{eq:Wstuv}
&&\sqrt{n}\Big| w_{st}^{\pm}-w_{uv}^{\pm} \Big| \\
&\leq& \Big| \cT^{\pm}(y;u,v, -\delta)-\cT^{\pm}(y;u,v, 0) \Big| + \Big| \cT^{\pm}(y;u,v, \delta)-\cT^{\pm}(y;u,v, 0) \Big|\nonumber\\
&\leq& \Big|\cZ^{\pm}(y;u,v, -\delta)-\cZ^{\pm}(y;u,v, 0) \Big|+ \Big|m^{\pm}(y;u,v, -\delta)- m^{\pm}(y;u,v, 0) \Big|\nonumber\\
&&+\Big|\cZ^{\pm}(y;u,v, \delta)-\cZ^{\pm}(y;u,v, 0) \Big|+ \Big|m^{\pm}(y;u,v, \delta)- m^{\pm}(y;u,v, 0) \Big|.\nonumber
\eenr

\begin{lem}\lel{lemma:31}
Assumption (B.2) implies
\ben
E\int\big[\cZ^{\pm}(y;u,v, \delta)-\cZ^{\pm}(y;u,v, 0) \big]^{2}d H(y)=o(1). \lel{eq:lemma311}
\een
Also, assumption (B.3) implies with $c$ and $\delta$ in (B.3)
\ben
\liminf_{n\ra\iny}P\Big( \sup \int n \big\{\nu^{\pm}(y,s,t)-\nu^{\pm}(y,u,v) \big\}^{2}dH(y)  \leq c\delta^{2}\Big)=1 \lel{eq:lemma312}
\een
where supremum is taken over $\|s-u\|\leq \delta$ and $\|t-v\|\leq \delta$. Moreover, assumptions (B.1)-(B.5) imply
\benr\lel{eq:lemma312_1}
&&\sup_{\|u\|\leq b,\|v\|\leq b}\int\,n^{-1}\Big[\sti  \vep_{i-k}\big\{ F_{\tiny{\xi}}(y+n^{-1/2}v'Z_{i}+\z_{ni}(u,v)) -F_{\tiny{\xi}}(y) \\
&&\quad\quad\quad\quad\quad\quad\quad\quad\quad\quad\quad  -f_{\tiny{\xi}}(y)\{ n^{-1/2}v'Z_{i}+\z_{ni}(u,v)\}\big\}  \Big]^{2}\,d H(y)=o_{p}(1).\nonumber
\eenr
\end{lem}
\noindent
\textbf{Proof}. Let $\cF_{i-1}:=\sigma\{\xi_{i-1},\xi_{i-2},...\}$. Define
\benn
\tilde{p}(y, u,v;z):= F_{\tiny{\xi}}(y+n^{-1/2}v'Z_{i}+  \z_{ni}(u,v)+d_{ni}(\delta,u,v)) - F_{\tiny{\xi}}(y+n^{-1/2}v'Z_{i}+  \z_{ni}(u,v)).
\eenn
Observe that $i$th summand is conditionally centered r.v.'s, given $\cF_{i-1}$, and hence, covariance of two summands is zero. Similar to (\r{eq:Bni}), the conditional variance of $i$th summand, given $\cF_{i-1}$, is bounded by
\benrr
E(\vep_{i-k}^{\pm})^{2}|\tilde{p}(y,u,v;Z_{i})| &\leq&  E(\vep_{i-k}^{\pm})^{2}| F_{\tiny{\xi}}(y+n^{-1/2}v'Z_{i}+  \z_{ni}(u,v)+d_{ni}(\delta,u,v)  ) - F_{\tiny{\xi}}(y)|\\
&&\quad\quad\quad\quad + E(\vep_{i-k}^{\pm})^{2}|F_{\tiny{\xi}}(y+n^{-1/2}v'Z_{i}+  \z_{ni}(u,v))- F_{\tiny{\xi}}(y) |.
\eenrr
Therefore, (\r{eq:lemma311}) follows from Fubini's Theorem and repeated application of (B.2) with $\delta$ and $\delta=0$. To prove (\r{eq:lemma312}), observe that the monotonicity of $F_{\tiny{\xi}}$ and (\r{eq:strange}) imply
\benr
&&m^{\pm}(y;u,v, -\delta)\leq m^{\pm}(y;s,t, 0)\leq m^{\pm}(y;u,v, \delta),\nonumber\\
&&m^{\pm}(y;u,v, -\delta)\leq m^{\pm}(y;u,v, 0)\leq m^{\pm}(y;u,v, \delta),\nonumber
\eenr
and hence,
\benn
m^{\pm}(y;u,v, -\delta)-m^{\pm}(y;u,v, 0) \leq n^{1/2}[\nu_{st}^{\pm}-\nu_{uv}^{\pm}]\leq m^{\pm}(y;u,v, \delta)-m^{\pm}(y;u,v, 0).
\eenn
Using the fact that $a_{1}\leq a_{2}\leq a_{3}$ implies $|a_{2}|\leq |a_{1}|+|a_{3}|$, we obtain
\benr\lel{eq:lemma313}
n^{1/2}|\nu_{st}^{\pm}-\nu_{uv}^{\pm}|&\leq& n^{-1/2}\sti  \vep_{i-k}^{\pm}\Big\{\,F_{\tiny{\xi}}(y+Z_{i}'v/\sqrt{n} +\z_{ni}(u,v)+d_{ni}(\delta, u,v)) \\
&&\qquad\qquad - F_{\tiny{\xi}}(y+Z_{i}'v/\sqrt{n} +\z_{ni}(u,v)-d_{ni}(\delta, u,v)) \,\Big\}.\nonumber
\eenr
As a consequence of (\r{eq:lemma313}) and the assumption (B.3), we finally obtain
\benr
&&\liminf_{n\ra\iny} P\Big( \sup \big| \sqrt{n}\{\nu^{\pm}(y,s,t)-\nu^{\pm}(y,u,v)\}\big|_{H}^{2}  \leq c\delta^{2}\Big)\nonumber\\
&\geq &\liminf_{n\ra\iny} P\Bigg( \sup \Big| n^{-1/2}\sti  \vep_{i-k}^{\pm}\Big\{\,F_{\tiny{\xi}}(y+Z_{i}'v/\sqrt{n} +\z_{ni}(u,v)+d_{ni}(\delta, u,v)) \nonumber\\
&&\quad\quad\quad - F_{\tiny{\xi}}(y+Z_{i}'v/\sqrt{n} +\z_{ni}(u,v)-d_{ni}(\delta, u,v)) \,\Big\}\Big|_{H}^{2}  \leq c\delta^{2}\Bigg)\nonumber\\
&=&1,\nonumber
\eenr
where the supremums are taken over $\|s-u\|\leq \delta$ and $\|t-v\|\leq \delta$, and hence, \eqref{eq:lemma312} follows. Next, write $m_{uv}^{\pm}$ and $m^{\pm}$ for $m^{\pm}(y;u,v,0)$ and $m^{\pm}(y;0,0,0)$. Let
\benn
\cL_{uv}=\big|m_{uv}^{\pm}-m^{\pm} - n^{-1/2}\sti  \vep_{i-k}^{\pm}f_{\tiny{\xi}}\{ n^{-1/2}v'Z_{i}+\z_{ni}(u,v)\}  \big|_{H}^{2}.
\eenn
In order to prove \eqref{eq:lemma312_1}, by using (B.4) and compactness of $\|u\|\leq b$ and $\|v\|\leq b$ (see Theorem 9.1.1 of Koul (2002)), it suffices to show that for all $\epsilon>0$ there exists a $\delta>0$ such that
\ben
\limsup_{n\ra\iny}P\Big(\sup\big|\cL_{st}-\cL_{uv}\big|>\epsilon \Big)=0,\lel{eq:lemma315}
\een
where sup is taken over $\|s\|,\|t\|\leq b$, $\|s-u \|\leq \delta$, and $||t-v ||\leq \delta$. Using monotonicity of $F_{\tiny{\xi}}$, \eqref{eq:strange}, triangle inequality, $(a_{1}\pm a_{2})^{2}\leq 2(a_{1}^{2}+a_{2}^{2})$, and (B.3), we obtain
for $\|s\|,\|t\|\leq b$ with $\|s-u \|\leq \delta$ and $\|t-v \|\leq \delta$
\ben\lel{eq:lemma317}
\limsup_{n\ra\iny} P\big(\sup|m_{st}-m_{uv}|_{H}^{2}\geq 4c\delta^{2}\big) = 0
\een
\noindent
where the equality follows from (\r{eq:lemma312}). Next, observe that
\benr
&&\quad\quad \big|\, \cL_{st}-\cL_{uv} \big|\nonumber\\
&\leq&  \int\,\Big[(m_{st}-m_{uv})- n^{-1/2}\sti  \vep_{i-k}\big\{ \big(Z_{i}'t/\sqrt{n} +\z_{ni}(s,t)\big) - \big(Z_{i}'v/\sqrt{n} +\z_{ni}(u,v)\big)\big\}f_{\tiny{\xi}}(y) \Big]\nonumber\\
&&\quad\quad\times \Big[ (m_{st}-m_{uv})+ 2\big\{ m_{uv}-m-n^{-1/2}\sti  \vep_{i-k}\big(Z_{i}'v/\sqrt{n} +\z_{ni}(u,v)\big) f_{\tiny{\xi}}(y) \big\}\nonumber\\
&&\quad\quad\quad - n^{-1/2}\sti  \vep_{i-k}\big\{ \big(Z_{i}'t/\sqrt{n} +\z_{ni}(s,t)\big) - \big(Z_{i}'v/\sqrt{n} +\z_{ni}(u,v)\big)\big\}f_{\tiny{\xi}}(y) \Big]d H(y)\nonumber
\eenr
\benr
&\leq& 2 \big|m_{st}-m_{uv}\big|_{H}^{2}+2\big|n^{-1/2}\sti  \vep_{i-k}d_{ni}(\delta,u,v)f_{\tiny{\xi}}\big|_{H}^{2}+2\big|m_{st}-m_{uv}\big|_{H}\nonumber\\
&&\times \big|n^{-1/2}\sti  \vep_{i-k}d_{ni}(\delta,u,v)\big|_{H} + 2\big|m_{uv}-m- n^{-1/2}\sti  \vep_{i-k}\big(Z_{i}'v/\sqrt{n} +\z_{ni}(u,v)\big) f_{\tiny{\xi}}\big|_{H}\nonumber\\
&&\quad\quad\quad\quad\quad\quad \times \Big\{ \big|m_{st}-m_{uv}\big|_{H}+\big|n^{-1/2}\sti  \vep_{i-k}d_{ni}(\delta,u,v)f_{\tiny{\xi}}\big|_{H}  \Big\}.\nonumber
\eenr
\noindent
In the first inequality, we use the fact that $(a_{1}^{2}-a_{2}^{2})=(a_{1}-a_{2})\cdot(a_{1}+a_{2})$; (\r{eq:lemma317}) implies that supremum of the first term of the last equation is $\delta o_{p}(1)$; the assumption (B.1) and the fact that $d_{ni}(\delta, u, v)=O_{p}(n^{-1/2}\delta\|Z_{i}\|)$ imply that supremum of  the second term is $\delta O_{p}(1)$, and hence, supremum of  the third term is $\delta o_{p}(1)$; the assumption (B.4) ensures that supremum of the last term is $\delta o_{p}(1)$. Therefore, by making $\delta$ as small as desired, (\r{eq:lemma315}) follows, thereby completing the proof of (\r{eq:lemma312_1}).

\begin{lem}\lel{lemma:31_1}
\benn
\sup_{\|u\|\leq b,\|v\|\leq b}\int\,\Big[ n^{-1/2}\sti  \vep_{i-k}\z_{ni}(u,v) f_{\tiny{\xi}}(y) \Big]^{2}\,d H(y)=o_{p}(1).\lel{eq:lemma31_11}
\eenn
\end{lem}
\noindent
\textbf{Proof}. Note that $\sup_{\|u\|\leq b,\|v\|\leq b}\z_{ni}(u,v)=O(n^{-1/2})$. Together with this fact, Fubini's theorem yields
\benr
&&\limsup_{n\ra\iny} E\sup_{\|u\|\leq b,\|v\|\leq b}\int n^{-1}\sti  \sum_{j=1}^{n}\vep_{i-k}\vep_{j-k}\z_{ni}(u,v)\z_{nj}(u,v)
f_{\tiny{\xi}}^{2}(y)\,d H(y)\nonumber\\
&\leq& \limsup_{n\ra\iny}\,C\,|f_{\tiny{\xi}}|_{H}^{2}\,\cdot n^{-2}\sti  \sum_{j=1}^{n}\alpha_{l}^{1/2}(j-i) \|\vep_{1}\|_{4}^{2} = 0,\nonumber
\eenr
where $0<C<\iny$. The inequality follows from Lemma \r{lemma:Deo}; the last equality follows from (B.1.a) and Lemma \r{lemma:Infinite_Test}. Consequently,
\benn
E\sup\big|  n^{-1/2}\sti  \vep_{i-k}\z_{ni}(u,v) f_{\tiny{\xi}}  \big|_{H}^{2} \ra 0,
\eenn
where the supremum is taken over $\|u\|\leq b$ and $\|v\|\leq b$, thereby completing the proof of lemma.
\begin{cor}\lel{cor:31}
\benn
\sup_{\|u\|\leq b,\|v\|\leq b }\big|m_{uv}-m_{0v}\big|_{H}^{2}=o_{p}(1).\lel{eq:cor311}
\eenn
\end{cor}
\noindent
\textbf{Proof}. Note that
\benrr
&& m(y,u,v,0 )-m(y,0,v,0 ) \\
&=&\Big\{  m(y,u,v,0 )-m(y,0,0,0 )-  n^{-1/2}\sti  \vep_{i-k} \big( Z_{i}'v/\sqrt{n} +\z_{ni}(u,v)\big)f_{\tiny{\xi}}(y)  \Big\}\\
&& - \Big\{  m(y,0,v,0 )-m(y,0,0,0 )-  v'n^{-1} \sti  \vep_{i-k}Z_{i}f_{\tiny{\xi}}(y)     \Big\} + n^{-1/2}\sti  \vep_{i-k} \z_{ni}(u,v)f_{\tiny{\xi}}(y).
\eenrr
Hence the proof of corollary will be completed by simply applying C-S inequality on the cross product terms after we show that
\benrr
&&\sup\,\big| m_{uv}-m-  n^{-1/2}\sti  \vep_{i-k} \big( Z_{i}'v/\sqrt{n} +\z_{ni}(u,v)\big)f_{\tiny{\xi}} \big|_{H}^{2}=o_{p}(1),\\
&&\sup\,\big| m_{0v}-m -  v' n^{-1} \sti  \vep_{i-k}Z_{i}f_{\tiny{\xi}}\big|_{H}^{2}=o_{p}(1),\\
&&\sup\, \big|  n^{-1/2}\sti  \vep_{i-k}\z_{ni}(u,v) f_{\tiny{\xi}} \big|_{H}^{2}=o_{p}(1),
\eenrr
where the supremums are taken over $\|u\|\leq b$ and $\|v\|\leq b$. The first and the second follow from (\r{eq:lemma312_1}) of Lemma \r{lemma:31}; the last follows from Lemma \r{lemma:31_1}.

\begin{lem}\lel{lemma:32}
\benn
\sup_{\|u\|\leq b,\|v\|\leq b }\int \Big\{ \cW^{\pm}(y,u,v)-\cW^{\pm}(y,0,0) \Big\}^{2}\,d H(y)=o_{p}(1).\lel{eq:lemma321}
\eenn
\end{lem}
\noindent
\textbf{Proof}. Similar to the square root of the integrand in (\r{eq:lemma311}), $\sqrt{n}\big[ \cW^{\pm}(y,u,v)-\cW^{\pm}(y,0,0) \big]$ is sum of conditionally centered r.v.'s. Thus, by repeated application of (B.2) with $\delta$ and $\delta=0$, we easily obtain for fixed $\|u\|,\,\|v\|\leq b$,
\ben
\big| \cW_{uv}^{\pm}-\cW^{\pm}\big|_{H}^{2}=o_{p}(1).\lel{eq:lemma321_1}
\een
Let $\cC_{uv}:=\big| \cW_{uv}^{\pm}-\cW^{\pm}\big|_{H}^{2}$. Similar to the proof of \eqref{eq:lemma315}, in order to complete the proof, it suffices to show that for every $\epsilon$ there exists $\delta$ such that for every $\|u\|,\,\|v\|\leq b$,
\ben
\limsup_{n\ra \iny} P\Big(\sup \big| \cC_{st} - \cC_{uv}\big|>\epsilon \Big)=0, \lel{eq:lemma322}
\een
where sup is taken over $\|s\|,\|t\|\leq b$, $\|s-u \|\leq \delta$, and $\|t-v \|\leq \delta$. Note that
\benn
\big|\, \cC_{st}-\cC_{uv} \big| \leq |\,\cW_{st}^{\pm}-\cW_{uv}^{\pm}\,|_{H}^{2}+2|\,\cW_{st}^{\pm}-\cW_{uv}^{\pm}|_{H}\cdot |\,\cW_{uv}^{\pm}-\cW^{\pm} \,|_{H}.
\eenn
By the monotonicity of $F_{\tiny{\xi}}$ together with \eqref{eq:strange} we obtain that
\benn
|m^{\pm}(y;u,v, \pm\delta)-m^{\pm}(y;u,v, 0)|\leq m^{\pm}(y;u,v, \delta) - m^{\pm}(y;u,v, -\delta).\lel{eq:lemma323}
\eenn
Recall definition of $\cW^{\pm}$. Then triangle inequality yields that
\benr
|\,\cW_{st}-\cW_{uv}\,|&\leq& |\,\cW_{st}^{+}-\cW_{uv}^{+}\,|+ |\,\cW_{st}^{-}-\cW_{uv}^{-}\,|,\nonumber\\
|\,\cW_{st}^{\pm}-\cW_{uv}^{\pm}\,|&\leq& \sqrt{n}\Big\{\,|w_{st}^{\pm}-w_{uv}^{\pm}|+ |\nu_{st}^{\pm}-\nu_{uv}^{\pm}| \Big\}.\nonumber
\eenr
Therefore, using $(a_{1}+a_{2})^{2}\leq 2(a_{1}^{2}+a_{2}^{2})$ repeatedly and \eqref{eq:Wstuv}, we have
\benr
|\,\cW_{st}^{\pm}-\cW_{uv}^{\pm}\,|_{H}^{2}&\leq& 16 \Bigg\{ \int[ \cZ^{\pm}(y;u,v, \delta)-\cZ^{\pm}(y;u,v, 0)  ]^{2}\,d H(y)\nonumber\\
&& +\int[ \cZ^{\pm}(y;u,v, -\delta)-\cZ^{\pm}(y;u,v, 0)  ]^{2}\,d H(y)\nonumber\\
&& +\int[ m^{\pm}(y;u,v, \delta)-m^{\pm}(y;u,v, -\delta)  ]^{2}\,d H(y)\nonumber\\
&& + \big|\sqrt{n}(\nu_{st}^{\pm}-\nu_{uv}^{\pm}) \big|_{H}^{2}\,\Bigg\}\nonumber
\eenr
\noindent
for all $\|s\|,\,\|t\|\leq $, $\|s-u \|\leq \delta$, $\|t-v \|\leq \delta$. Finally, assumption (B.3), \eqref{eq:lemma311}, \eqref{eq:lemma312},  and \eqref{eq:lemma321_1} prove \eqref{eq:lemma322}, thereby completing the proof of the lemma.
\vs .15cm
Define
\benr
SS_{k}(y,u,v)&:=& n^{-1/2}\sti   \vep_{i-k} \Big\{\,I\big( \xi_{i} \leq y+Z_{i}'v/\sqrt{n} +\z_{ni}(u,v)\big)\nonumber\\
 &&\quad\quad -I\big( -\xi_{i}< y-Z_{i}'v/\sqrt{n} -\z_{i}(u,v)\big)\,\Big\}.\nonumber
\eenr
\noindent
Rewrite
\benr
\,\,\,\,R1_{uv}^{k}&=& \Big[\cW(y,u,v)-\cW(y,0,0)\Big] - \Big[\cW(y,0,v)-\cW(y,0,0)\Big]\nonumber\\
&&+ \Big[\cW(-y,u,v)-\cW(-y,0,0)\Big] - \Big[\cW(-y,0,v)-\cW(-y,0,0)\Big]\nonumber\\
&& + \Big[ m(y,u,v,0)-m(y,0,v,0)\Big] - \Big[ m(-y,u,v,0)-m(-y,0,v,0)\Big].\nonumber
\eenr
Using quadratic expansion and applying C-S inequality involving $L_{H}^{2}$ norm on the cross product, the claim that $\sup_{\|u\|,\|v\|}R1_{uv}^{k}=o_{p}(1)$ will be true if we can show that
\benrr
&&\textrm{(a)}\sup_{\|u\|\leq b,\|v\|\leq b}|\cW_{uv}- \cW |_{H}^{2}=o_{p}(1),\quad \textrm{(b)}\sup_{\|v\|\leq b}|\cW_{0v}- \cW |_{H}^{2}=o_{p}(1),\\
&&\textrm{(c)}\sup_{\|u\|\leq b,\|v\|\leq b}\left| m_{uv}-m_{0v}\right|_{H}^{2}=o_{p}(1).
\eenrr
(a) and (b) follow from Lemma \r{lemma:32}; (c) follows from Corollary \r{cor:31}.

Next, consider $R2_{uv}^{k}$. To begin with, note that by C-S inequality and Fubini's theorem together, we obtain for $v\in\R^{q}$, and $0<b<\iny$
\benr\lel{eq:R00_1}
&&E\int\Big\{F_{\tiny{\xi}}(y+v'Z_{1}/\sqrt{n})-F_{\tiny{\xi}}(y-v'Z_{1}/\sqrt{n})\Big\}^{2}d H(y)\\
&\leq&  4b^{2} n^{-1} \big(2b/\sqrt{n}\big)^{-1}\int_{-b/\sqrt{n}}^{b/\sqrt{n}}\int E\big[\|Z_{1}\|f_{\tiny{\xi}}(y+s\|Z_{1}\|)\big]^{2}\,d H(y)\,ds.\nonumber
\eenr
Moreover, by the assumptions (B.1.a) and (B.1.c) and bounded convergence theorem, we have
\benr \lel{eq:R00_2}
&&\int\Big[E\big|F_{\tiny{\xi}}(y+ v'Z_{1}/\sqrt{n})-F_{\tiny{\xi}}(y-v'Z_{1}/\sqrt{n})\big| \Big]^{1/2}d H(y)\\
&\leq&\sqrt{2b}\,n^{-1/4}\int\Big[ (2b/\sqrt{n})^{-1}\int_{-b/\sqrt{n}}^{b/\sqrt{n}}E\{\|Z_{1}\|f_{\tiny{\xi}}(y+s\|Z_{1}\|)\}ds\Big]^{1/2}d H(y) =O(n^{-1/4}).\nonumber
\eenr
Finally, consider $v_{1},\,v_{2}\in \R^{q}$ and let $\|v_{1}\|\leq \|v_{2}\|$. Using C-S inequality and Fubini's theorem again, we have
\benr \lel{eq:R00_3}
&&\,\,\,\, E\int\Big\{F_{\tiny{\xi}}(y+v_{2}'Z_{1}/\sqrt{n})-F_{\tiny{\xi}}(y+v_{1}'Z_{1}/\sqrt{n})\Big\}^{2}d H(y)\\
&\leq&  n^{-1} ( \|v_{2}\|+\|v_{1}\|)^{2}  \big\{(\|v_{2}\|+\|v_{1}\|)/\sqrt{n}\big\}^{-1}\int_{-\|v_{1}\|/\sqrt{n}}^{\|v_{2}\|/\sqrt{n}}\int E\big[\|Z_{1}\|f_{\tiny{\xi}}(y+s\|Z_{1}\|)\big]^{2}\,d H(y)\,ds.\nonumber
\eenr
Define
\benn
\mbf{q}(y,v;z):=F_{\tiny{\xi}}(y+n^{-1/2}v'z)-F_{\tiny{\xi}}(y-n^{-1/2}v'z).
\eenn
\noindent
Next, rewrite
\benr
R2_{uv}^{k}&=& n^{-1/2} \sti   \eta_{u}^{ik} \Big\{\,I\big(\xi_{i}  \leq y+n^{-1/2}v'Z_{i} \big) - I\big(-\xi_{i} < y-n^{-1/2}v'Z_{i} \big) - \mbf{q}(y,v;Z_{i}) \,\Big\}\nonumber\\
&&  + n^{-1/2} \sti   \eta_{u}^{ik}\,\mbf{q}(y,v;Z_{i}) \nonumber\\
&=& R21_{uv}^{k}+R22_{uv}^{k}, \quad\quad say.\nonumber
\eenr
Due to the monotonicity of $F_{\tiny{\xi}}$, we obtain
\ben
\sup_{\|v\|\leq b}|\mbf{q}(y,v;Z_{i})|\leq F_{\tiny{\xi}}(y+n^{-1/2}b\|Z_{i}\|)-F_{\tiny{\xi}}(y-n^{-1/2}b\|Z_{i}\|).\lel{eq:R0_1}
\een
Therefore, by C-S inequality, (\r{eq:R0_1}), $\va_{n}=o(1)$, and (\r{eq:R00_1}), we have
\benn
E\,\sup_{\|u\|,\|v\|\leq b }|R22_{uv}^{k}|_{H}^{2} \leq  4b^{4}\,\va_{n}^{2}\,\big(2b/\sqrt{n}\big)^{-1}\int_{-b/\sqrt{n}}^{b/\sqrt{n}}\int E\big[\|Z_{1}\|f_{\tiny{\xi}}(y+s\|Z_{1}\|)\big]^{2}\,d H(y)\,ds \ra 0.
\eenn
Next, consider the cross product of $R21_{uv}^{k}$. For fixed $u,\,v$
\benr
&& E\,\int \frac{1}{n}\sti  \sum_{j=i+1}^{n} \Big[\, \eta_{u}^{ik}\eta_{u}^{jk}\Big\{\,I\big(\xi_{i}  \leq y+n^{-1/2}v'Z_{i} \big) - I\big(-\xi_{i} < y-n^{-1/2}v'Z_{i} \big)  \nonumber\\
&& \quad - \mbf{q}(y,v;Z_{i}) \,\Big\} \times \Big\{\,I\big(\xi_{j}  \leq y+n^{-1/2}v'Z_{j} \big) - I\big(-\xi_{j} < y-n^{-1/2}v'Z_{j} \big)\nonumber\\
&&\quad\quad \quad\quad - \mbf{q}(y,v;Z_{j}) \,\Big\} \Big]\, d H(y)\nonumber\\
&\leq& b^{2}\va_{n}^{2}{n}^{-1}\sti  \sum_{j=i+1}^{n}\int \Big|E\Big\{\,I\big(\xi_{i}  \leq y+n^{-1/2}v'Z_{i} \big) - I\big(-\xi_{i} < y-n^{-1/2}v'Z_{i} \big)\nonumber  \\
&& \quad -\mbf{q}(y,v;Z_{i}) \,\Big\} \times \Big\{\,I\big(\xi_{j}  \leq y+n^{-1/2}v'Z_{j} \big) - I\big(-\xi_{j} < y-n^{-1/2}v'Z_{j} \big)\nonumber\\
&&\quad \quad - \mbf{q}(y,v;Z_{j}) \,\Big\} \Big|d H(y)\nonumber\\
&\leq& 10b^{2}\,\va_{n}^{2}\,n^{-1}\sti  \sum_{j=i+1}^{n}\alpha_{l}^{1/2}(j-i)\nonumber\\
&& \quad\quad \quad\quad \times \int\Big[E\big|F_{\tiny{\xi}}(y+ v'Z_{1}/\sqrt{n})-F_{\tiny{\xi}}(y-v'Z_{1}/\sqrt{n})\big| \Big]^{1/2}d H(y)\nonumber\\
&\ra& 0.\nonumber
\eenr
The first inequality follows from Fubini's theorem; the second inequality follows from Lemma \r{lemma:Deo}; the convergence to zero follows from $\va_{n}=0(1)$, \eqref{eq:lemma_InfTest} with $r=2$, and \eqref{eq:R00_2}. Therefore, the sum of the covariance of two summands in $R21_{uv}^{k}$ will tend to zero. Consequently, for fixed $u$ and $v$, by Fubini's theorem, stationarity of $\vep$, and (B.5), we obtain
\benr\lel{eq:R21_2}
&&\limsup_{n\ra\iny}E|R21_{uv}^{k}|_{H}^{2}\\
&\leq&\limsup_{n\ra\iny} 2b^{3}\,\va_{n}^{2}\,n^{-1/2}\, (2b/n^{1/2})^{-1}\int_{-b/n^{1/2}}^{b/n^{1/2}}\int E\|Z_{1}\|f_{\tiny{\xi}}(y+s\cdot \|Z_{1}\|)d H(y)ds\nonumber\\
&=&0.\nonumber
\eenr

Finally, we shall show that $E\sup |R21_{uv}^{k}|_{H}^{2}=o(1)$. Let $u_{h}=x_{h}\,e_{u}$ and $v_{l}=w_{l}\, e_{v}$ for $0\leq l,\,h\leq r$ where $x_{h}=\|u_{h}\|$ and $w_{l}=\|v_{l}\|$, and $ e_{u},\,  e_{v}$ are unit vectors in $\R^{p+1}$ and $\R^{q}$. Let
\benn
-b=x_{0}\leq x_{1}\leq\cdots\leq x_{r}=b,\quad  -b=w_{0}\leq w_{1}\leq\cdots\leq w_{r}=b
\eenn
so that $\max_{h}(x_{h}-x_{h-1})\ra 0$ and $\max_{l}(w_{l}-w_{l-1})\ra 0$ as $r\ra \iny$. Note that $\eta_{u}^{ik}=u'c_{i-k}$ and rewrite
\benr
R21_{uv}^{k}&=& n^{-1/2} \sti   u'c_{i-k} \Big\{\,I\big(\xi_{i}  \leq y+n^{-1/2}v'Z_{i} \big)- I\big(-\xi_{i} < y-n^{-1/2}v'Z_{i} \big) \,\Big\}\nonumber\\
&&\quad\quad - n^{-1/2} \sti   u'c_{i-k}\,\mbf{q}(y,v;Z_{i}) \nonumber\\
&=& T1_{uv}-T2_{uv},\quad say.\nonumber
\eenr
Define $\eta_{u}^{ik+}:=\eta_{u}^{ik}I(\eta_{u}^{ik}Z_{i}'v>0)=u'c_{i-k}I(u'c_{i-k}Z_{i}'v>0)$, $\eta_{u}^{ik-}=\eta_{u}^{ik}-\eta_{u}^{ik+}$. Let $R21_{h,l}^{\pm}$, $T1_{h,l}^{\pm}$, $T2_{h,l}^{\pm}$ stand for the $R21_{u_{h}v_{l}}^{k}$, $T1_{u_{h}v_{l}}$, $T2_{u_{h}v_{l}}$ with $\eta_{u}^{ik}$ being replaced by $\eta_{u}^{ik\pm}$. Observe that
\benn
|R21_{uv}^{k}|_{H}^{2}\leq 2\Big\{ |T1_{uv}^{+}-T2_{uv}^{+}|_{H}^{2}+|T1_{uv}^{-}-T2_{uv}^{-}|_{H}^{2} \Big\}.
\eenn
Note that $R21_{uv}^{+}\,(R21_{uv}^{-})$ is a difference of two nondecreasing (nonincreasing) function of $x,\,w$, and hence, for all $u=x\,e_{u}\in\R^{p+1}$ and $v=w\,e_{v}\in\R^{q}$ where $x_{h-1}\leq x\leq x_{h}$ and $w_{l-1}\leq w\leq w_{l}$,
\benn
R21_{h-1,l-1}^{+}-(T2_{h,l}^{+}-T2_{h-1,l-1}^{+})\leq R21_{uv}^{+}\leq R21_{h,l}^{+}+(T2_{h,l}^{+}-T2_{h-1,l-1}^{+}).
\eenn
Therefore, using the fact that $a_{1}\leq a_{2}\leq a_{3}$ implies $|a_{2}|\leq |a_{1}|+|a_{3}|$ and that $(a\pm b)^{2}\leq 2(a^{2}+b^{2})$ repeatedly, we obtain
\benr \lel{eq:supR21}
\sup |R21_{uv}^{+}|_{H}^{2} &\leq& 8\Big[\max_{0\leq h,l\leq r}|R21_{h,l}^{+}|_{H}^{2}+ 2\max_{0\leq h,l\leq r}|T2_{h,l}^{+}-T2_{h-1,l}^{+}|_{H}^{2}\\
&&\quad\quad  + 2\max_{0\leq h,l\leq r}|T2_{h-1,l}^{+}-T2_{h-1,l-1}^{+}|_{H}^{2}  \Big]\nonumber.
\eenr
Using $\|\eta_{u}^{ik+}\|^{2}\leq \|\eta_{u}^{ik}\|^{2}$ and ergodicity of $\vep$,
\benr\lel{eq:T2pm1}
&&E|T2_{h,l-1}^{+}-T2_{h-1,l-1}^{+}|_{H}^{2}\\
&\leq& 16b^{4}\,\va_{n}^2\, \big(2b/\sqrt{n}\big)^{-1}\int_{-b/\sqrt{n}}^{b/\sqrt{n}}\int E\big[\|Z_{1}\|f_{\tiny{\xi}}(y+s\|Z_{1}\|)\big]^{2}\,d H(y)\,ds\ra 0,\nonumber
\eenr
where the inequality follows from (\r{eq:R00_1}), and the convergence to zero follows from assumption (B.5) and $\va_{n}=o(1)$. Similarly, using (\r{eq:R00_3}), we have
\benr \lel{eq:T2pm2}
&&E|T2_{h-1,l}^{+}-T2_{h-1,l-1}^{+}|_{H}^{2}\\
 &\leq& b^{2}\va_{n}^2 (|w_{l}|+|w_{l-1}|)^{2}\cdot \big\{(|w_{l}|+|w_{l-1}|)/\sqrt{n}\big\}^{-1} \nonumber \\
&&\times \Bigg\{ \int_{-|w_{l-1}|/\sqrt{n}}^{|w_{l}|/\sqrt{n}}\int E\big[\|Z_{1}\|f_{\tiny{\xi}}(y+s\|Z_{1}\|)\big]^{2}\,d H(y)\,ds \nonumber\\
 &&\quad\quad\quad\quad +\int_{-|w_{l-1}|/\sqrt{n}}^{|w_{l}|/\sqrt{n}}\int E\big[\|Z_{1}\|f_{\tiny{\xi}}(y-s\|Z_{1}\|)\big]^{2}\,d H(y)\,ds\Bigg\}\nonumber\\
&\ra& 0.\nonumber
\eenr
Consequently, \eqref{eq:R21_2}, \eqref{eq:supR21}, \eqref{eq:T2pm1}, and \eqref{eq:T2pm2} imply
\benn
E\sup_{u,v} |R21_{uv}^{+}|_{H}^{2}\longrightarrow 0.
\eenn
Similar facts hold for $|R21_{uv}^{-}|_{H}^{2}$, and this completes the proof of $\sup R2_{uv}^{k}=o_{p}(1)$.

To show $\sup_{\|u\|\leq b,\|v\|\leq b} \big|R3_{uv}^{k}\big|_{H}^{2}$ is $o_{p}(1)$, note that $R3_{uv}^{k}$ is the same as $R1_{uv}^{k}$ except the fact that it has $o(1)$ as a weight instead of random error. So, the proof is similar, but much simpler. Hence we do not present here.

\end{document}